\documentclass[preprint,12pt]{elsarticle}

\usepackage{amssymb}
\usepackage{amsmath}
\usepackage{graphicx}
\usepackage{subcaption}
\usepackage{float,ulem}
\usepackage{mathtools}
\usepackage{url,hyperref}
\usepackage{booktabs}
\usepackage{xcolor,todonotes}
\usepackage{placeins}

\makeatletter
\AtBeginDocument{\let\old@subsection\subsection
  \renewcommand{\subsection}{\FloatBarrier\old@subsection}}
\makeatother

\definecolor{zzrev}{RGB}{0,0,0}
\newcommand{\revzz}[1]{{\color{zzrev}#1}}
\newcommand{\revaa}[1]{{\color{zzrev}#1}}

\journal{Computer Methods in Applied Mechanics and Engineering}

\begin{document}

\begin{frontmatter}



\title{
Adjoint Method versus Physics-Informed Neural Networks in PDE-Constrained Inverse Problems
}


\author[label1]{Zhen Zhang}
\ead{zhen\_zhang2@brown.edu}
\author[label2]{Alessandro Alla}
\ead{alessandro.alla@uniroma1.it}
\author[label1]{George Em Karniadakis\corref{cor1}}
\ead{george_karniadakis@brown.edu}
\cortext[cor1]{Corresponding author}

\affiliation[label1]{organization={Division of Applied Mathematics, Brown University},
            city={Providence},
            state={RI},
            country={USA}}

\affiliation[label2]{organization={Department of Mathematics, Sapienza Università di Roma},
            city={Rome},
            country={Italy}}

\begin{abstract}

Inverse problems governed by partial differential equations (PDEs) are central
to computational mechanics and are commonly solved by adjoint-based
optimization, while physics-informed neural networks (PINNs) have emerged as a
flexible alternative. Their relative performance remains difficult to assess
because the two approaches are often compared under different formulations,
parameterizations, optimizers, and regularization choices.
We present a \revzz{controlled} comparison of adjoint optimization and PINNs for
PDE-constrained inverse problems. From a common abstract formulation, we
instantiate both methods on identical domains, governing equations,
observation models, and regularization terms, while matching the optimizer,
unknown parameterization, and arithmetic precision wherever applicable. The
benchmarks include unsteady Burgers, noisy Darcy permeability inversion,
three-dimensional Allen--Cahn reaction identification, and unsteady
Navier--Stokes viscosity identification.
The results show that the representation of the unknown largely determines
the preferred method: grid-based fields favor the discrete adjoint, whereas
neural representations are native to PINNs and relevant for closure and
constitutive modeling. For time-dependent problems, adjoint inversion can be
dominated by trajectory storage and differentiation, while PINNs provide
satisfactory reconstructions at lower cost. A PINN-warm-started adjoint
strategy then recovers adjoint-level accuracy at \revzz{about half the cost}.

\end{abstract}


\begin{highlights}
\item Rigorous adjoint and PINN inversion compared on PDE benchmarks
\item Grid representations favor the adjoint; neural networks favor PINNs
\item In unsteady problems, adjoint trajectory storage dominates the cost
\item PINN warm starts polished by adjoint reach target accuracy at under half the cost
\end{highlights}


\end{frontmatter}

\section{Introduction}

Inverse problems governed by partial differential equations (PDEs) are central
to computational mechanics, with applications ranging from parameter
identification and source reconstruction to model calibration
\cite{lions1971optimal,hinze2009optimization}. Traditionally, the unknown to
be inferred is a scalar parameter or a spatially distributed field. In many
emerging applications, however, the target is a more complex object: hidden physical parameterizations or state-dependent functions inferred directly from data.
This shift is particularly
prominent in data-driven closure and constitutive modeling across fluid and
solid mechanics \cite{duraisamy2019turbulence,karniadakis2021physics}, where
the unknown is often high-dimensional, state-dependent, and more naturally
represented by a neural network than by a grid-based field. This raises a
basic methodological question: does an inverse solver remain efficient and
accurate when the unknown moves from a discretized field to a neural
functional representation? This question motivates the present study.\\
\revaa{At the continuous level, these inverse problems are naturally posed in infinite-dimensional function spaces, with the unknowns  belonging, for instance, to a suitable space of spatial fields or nonlinear functions. Computational inversion therefore requires a finite-dimensional representation of this unknown. Grid-based discretizations, truncated basis expansions such as Karhunen–Loève representations, and neural networks can all be interpreted as different finite-dimensional approximations of the underlying infinite-dimensional inverse problem. The present work focuses on the computational consequences of this choice of representation, rather than on convergence with respect to the representation dimension.}
Adjoint-based gradient optimization has long been the standard approach for
PDE-constrained inverse problems
\cite{lions1971optimal,hinze2009optimization,jameson1988aerodynamic,giles2000introduction}.
In a discretize-then-optimize formulation, the adjoint of the fully discrete
forward problem provides the reduced gradient at machine precision through a
single backward sweep, with a cost essentially independent of the number of
optimization variables. This efficiency at the gradient level is the main
strength of the method, but the inverse problem is solved only after an outer
direct-adjoint loop (DAL): a forward solve, an adjoint solve, and a parameter update
are repeated until convergence
\cite{skene2021parallel,nzoyem2023comparison}. The total cost therefore
depends not only on the cost of one forward--adjoint sweep, but also on the
number of optimization iterations. Moreover, this efficiency comes with an
implementation cost: for each forward discretization, the corresponding
discrete adjoint must be derived, implemented, and verified. In particular,
the adjoint solver requires the transpose of the linearized discrete forward
operator and is tied to the chosen time integrator, spatial discretization,
and treatment of boundary conditions. Closely related to identification and control in PDE-constrained optimization, including also feedback control was studied in e.g. ~\cite{APPP,AP25}.

Physics-informed neural networks (PINNs) have recently emerged as an
alternative framework for forward and inverse problems governed by
PDEs~\cite{raissi2019physics,karniadakis2021physics}. Instead of solving the
PDE repeatedly on a fixed mesh, a PINN represents the state with one neural
network and the unknown with another, and enforces the governing equations
through a soft penalty on the PDE residual evaluated by automatic
differentiation~\cite{baydin2018automatic}. The state and the unknown are
then inferred jointly from the physics residual, boundary and initial
conditions, and observation data. This mesh-free formulation requires only the
PDE residual, rather than a discretization-specific adjoint or hand-derived
linearization, and is therefore naturally suited to unknowns represented by
neural networks, such as data-driven closures or constitutive laws. It also
replaces the outer direct-adjoint loop by a single joint optimization problem.
These advantages come with weaker accuracy and convergence guarantees, and
with training difficulties associated with ill-conditioned PINN losses.
Recent second-order quasi-Newton optimizers, such as
SSBroyden~\cite{urban2025unveiling,KIYANI2025118308,jnini2026curvature,bioli2026selfscaled}, have narrowed this gap. However, the
regimes in which PINNs provide a genuine practical advantage over adjoint
optimization remain incompletely characterized.

Despite a rapidly expanding literature on PINNs, adjoint methods, and
PDE-constrained optimization, a systematic comparison of these paradigms for
inverse discovery problems remains incomplete. Early PINN-based approaches
addressed constant-parameter identification and PDE control problems
\cite{alla2025pinn}, and recent comparisons between PINNs and adjoint-based
methods have clarified important aspects of their relative performance
\cite{nzoyem2023comparison,mowlavi2023optimal}. However, several issues remain
open.  
First, many existing studies are framed primarily as optimal control problems,
where the governing PDE is fully known and the objective is to determine a
control that minimizes a cost. This differs from inverse discovery problems,
where the unknown enters the operator itself, for example as a permeability, a
state-dependent reaction term, or a physical coefficient. Second, most
comparisons fix a single representation of the unknown, making it difficult to
separate the effect of the representation from that of the inversion algorithm.
Third, algorithmic differences are often entangled with implementation choices,
such as the optimizer, arithmetic precision, or parameterization used by the two
methods. Finally, time-dependent and higher-dimensional regimes, where the
storage and differentiation of the full state trajectory can dominate adjoint
cost, have received comparatively less attention. As a result, the literature
shows that performance is problem-dependent, but does not fully identify which
problem properties determine the preferred method.

The central premise of this work is that the preferred method can often be
anticipated from two properties of the inverse problem. The first is the
representation and structure of the unknown. A grid-defined or coarse-grid
field is naturally aligned with a discrete adjoint solver, whereas a neural
representation is native to PINNs and is the natural setting for discovering
data-driven closures or constitutive laws. In an adjoint formulation, the same
neural representation moves the optimization into a high-dimensional weight
space, where quasi-Newton steps can become expensive and poorly conditioned.
The second property is the cost of differentiating time-dependent trajectories.
For unsteady PDEs, each adjoint iteration requires a forward solve, a backward
adjoint sweep, and access to the full state history. This cost grows with the
spatial resolution, the number of time steps, and the dimension of the
discretized state, and can become the dominant bottleneck. PINNs avoid this
discrete trajectory bottleneck, although their cost is then governed by the
network architecture, collocation budget, and optimization dynamics.

Building on these observations, we present a \revzz{controlled} comparison of adjoint optimization and PINNs for four
benchmark PDE-constrained inverse problems of increasing nonlinearity and
dimensionality: forcing identification in the unsteady one-dimensional viscous
Burgers equation, log-permeability identification in steady two-dimensional
Darcy flow with sparse noisy observations, state-dependent reaction
identification in the unsteady three-dimensional Allen--Cahn equation, and
scalar-viscosity identification in the unsteady two-dimen\-sional
Navier--Stokes flow past a cylinder. 
All benchmarks are derived from the same abstract inverse-problem formulation.
For each case, the domain, governing equation, observation operator, and
regularization are kept identical across the adjoint and PINN formulations.
Where applicable, we also match the SSBroyden optimizer, the parameterization of the unknown, and the arithmetic precision.

\revzz{One difference cannot be matched in the comparison. The two methods do not solve the same
finite-dimensional problem. The adjoint imposes the discrete forward map as a hard
constraint, optimizing on the manifold where the discrete residual vanishes
exactly; the PINN imposes it softly, as a penalty in a loss minimized over the
state and parameter networks jointly, so the governing equation holds only
approximately at the optimum. Both are legitimate finite-dimensional realizations
of the same continuous inverse problem, but no choice of settings makes them
identical, and removing that difference would remove exactly what the study
measures. Therefore, we control the comparison in the sense of matching every factor that does admit a common
setting so that the differences that remain are attributable to the two paradigms rather than to incidental implementation choices.}

The experiments lead to three main findings. First, across the benchmarks
considered here, the representation of the unknown---as a grid-based field or
as a neural network---is a primary determinant of relative performance, and is
directly connected to the setting of data-driven closure and constitutive
modeling. Second, in time-dependent problems, the storage and differentiation
of the discrete state trajectory can make adjoint-based inversion expensive;
in these regimes, PINNs can recover the unknown to satisfactory accuracy at
substantially lower wall-clock cost. Third, when the accuracy of the discrete
adjoint is desired, a hybrid strategy in which the PINN solution is used to
warm start the adjoint optimization can recover high-fidelity accuracy at a
fraction of the cost of a cold-started adjoint solve.

This last construction---a lower-cost, lower-fidelity PINN reconstruction
polished by the discrete adjoint---is consistent with the broader philosophy
of multilevel optimization~\cite{gratton2008recursive} and multiscale
continuation for inverse problems~\cite{bunks1995multiscale}, and has recently
been advocated in PINN-based optimal control~\cite{zhang2026pinns}. Here, we
bring this idea to PDE-constrained inverse discovery in a controlled,
gradient-verified setting.

The remainder of the paper is organized as follows. Section~\ref{sec:inverse}
states the abstract PDE-constrained inverse problem and describes the discrete
adjoint and PINN paradigms in a unified notation. Section~\ref{sec:benchmarks}
instantiates it on the four benchmarks and reports the recovered fields,
accuracies, and wall-clock costs. Conclusions and practical guidance are drawn
in Section~\ref{sec:conclusion}. \revzz{\ref{sec:appendix} collects the
implementation details --- discretization, parameterization, regularization and
optimizer settings --- for each benchmark, and~\ref{sec:ablation} reports the
ablation studies, covering the choice of optimizer and of network architecture.}
\section{Inverse Problems and Solution Methods}
\label{sec:inverse}
 
This section formulates the class of PDE-constrained inverse problems studied in this work and describes the two solution paradigms whose performance is compared throughout the paper.
We first state a generic abstract formulation that encompasses all four benchmark problems considered, later, in Section~\ref{sec:benchmarks}.
We then describe the discrete adjoint approach and the PINN approach in a unified notation, emphasizing the design choices that are kept identical across the two methods \revzz{so that the numerical comparisons are free of the implementation artifacts that would otherwise confound them, and the differences that remain are properties of the two paradigms}.
 
\subsection{General Inverse Problem Formulation}
\label{subsec:general_inverse}
 
Let $\Omega\subset\mathbb{R}^d$ be a bounded spatial domain and $[0,T]$ a time interval. 
We consider PDE systems of the abstract form 
\begin{align}
\begin{aligned}
u_t+\mathcal{N}\bigl(u;f\bigr)&=0 \quad\text{in }\Omega\times(0,T], 
\\
\mathcal{B}(u)&=0 \quad\text{on }\partial\Omega\times(0,T],\\
u(\cdot,0)&=u_0 \quad\text{in }\Omega,
\label{eq:abstract_pde}
\end{aligned}
\end{align}
where $u\in\mathcal{X}$ denotes the state variable in the function space $\mathcal{X} = L^2(\Omega \times [0,T])$, $f$ is the unknown function in PDE to be discovered, and $\mathcal{N},\mathcal{B}$ are, respectively, differential and boundary operators encoding the governing physics.
Equation~\eqref{eq:abstract_pde} defines the \emph{forward map} $f\mapsto u(f)$, which we assume to be well-posed.

The inverse problem consists of recovering $f$ from observations.
We denote by $\mathcal{O}: \mathcal{X} \to \mathcal{Y}$ an observation operator mapping from the state space $\mathcal{X}$ to the observation space $\mathcal{Y}$, so that the available data are
\begin{equation}
y=\mathcal{O}\bigl(u(f^\star)\bigr)+\eta,
\label{eq:obs_model}
\end{equation}
where $f^\star$ is the unknown ground truth and $\eta$ models measurement noise.
The inverse problem is formulated as the regularized minimization
\begin{equation}
\min_{f}\; J(f)=\frac{1}{2}\,\bigl\|\mathcal{O}\bigl(u(f)\bigr)-y\bigr\|_{\mathcal{Y}}^{2} +\mathcal{R}(f), \quad \text{subject to } u(f) \text{ solves } \eqref{eq:abstract_pde},
\label{eq:abstract_inverse}
\end{equation}
where $\|\cdot\|_{\mathcal{Y}}$ denotes the norm in the observation space, and $\mathcal{R}(f)$ is a regularization functional.
 
Two distinct elements of \eqref{eq:abstract_inverse} drive the methodological discussion of this paper.
First, the unknown function $f$ is approximated via a finite-dimensional parameterization, denoted as $f(\boldsymbol{\theta}_f)$, where $\boldsymbol{\theta}_f \in \mathbb{R}^{N_f}$ represents the underlying degrees of freedom and $N_f$ denotes the total number of tunable parameters. The choice of parameterization controls the effective dimensionality of the problem and provides a first source of implicit regularization. In this work, $\boldsymbol{\theta}_f$ may be structured in three distinct ways: (i) as a grid-defined vector consistent with the forward discretization, (ii) as a coarse-grid-defined vector interpolated to the forward simulation mesh via a prolongation operator, or (iii) as the output of an over-parameterized neural network where $\boldsymbol{\theta}_f$ denotes the network weights.

Second, two algorithmic paradigms can be used to solve \eqref{eq:abstract_inverse}: \emph{adjoint-based optimization} \cite{lions1971optimal,jameson1988aerodynamic}, which iterates on the discrete forward solver and provides exact gradients of $J$ via a backward-in-time sweep, and \emph{physics-informed neural networks} \cite{raissi2019physics}, which replace the forward solve by a neural surrogate and minimize a composite physics-informed loss.

\subsection{Adjoint-Based Optimization}
\label{subsec:adjoint}
In this section, we recall the adjoint-based optimization. The interested reader can find an extensive description in e.g. \cite{hinze2009optimization,giles2000introduction}.
The adjoint approach treats \eqref{eq:abstract_inverse} as a PDE-constrained optimization problem and computes the reduced gradient $\nabla_{\boldsymbol{\theta}_f} J\bigl(f(\boldsymbol{\theta}_f)\bigr)$ at machine precision via a single backward sweep through a discretized adjoint equation.
This subsection summarizes the discrete forward problem, the discrete adjoint formulation, the verification protocol, and the optimizer.

\paragraph{Discrete forward problem}
Let $\mathbf{u}=(\mathbf{u}^0,\mathbf{u}^1,\dots,\mathbf{u}^{N_t})$ denote the discrete state trajectory produced by a chosen forward solver applied to \eqref{eq:abstract_pde}, where superscript $n$ denotes the $n$-th time step such that $\mathbf{u}^n \in \mathbb{R}^{N_x}$ for all $n \in \{0, 1, \dots, N_t\}$.

We work exclusively with the \emph{discretize-then-optimize} paradigm: the discrete adjoint is derived by differentiating the fully discrete forward residual vector $\mathbf{R}^n$. By replacing the temporal derivative $u_t$ and the spatial differential operator $\mathcal{N}(u; f)$ in \eqref{eq:abstract_pde} with their respective discrete counter-parts, the algebraic residual vector $\mathbf{R}^n: \mathbb{R}^{N_x} \times \mathbb{R}^{N_x} \times \mathbb{R}^{N_f} \to \mathbb{R}^{N_x}$ at the $n$-th time step is defined as:
\begin{equation}
\mathbf{R}^n\bigl(\mathbf{u}^n,\mathbf{u}^{n-1};\boldsymbol{\theta}_f\bigr) = \mathbf{0}, \qquad n=1,2,\dots,N_t,
\label{eq:discrete_residual}
\end{equation}
where $N_x$ denotes the number of spatial degrees of freedom resulting from the discretization mesh and $N_t$ denotes the number of time steps. This algebraic formulation implicitly accounts for both the internal physics and the boundary conditions $\mathcal{B}(u)=0$. Consequently, the continuous optimization problem \eqref{eq:abstract_inverse} reduces to a finite-dimensional, constrained optimization problem. 
By introducing the parameterization $f(\boldsymbol{\theta}_f)$, the target is to minimize the discrete objective function directly with respect to the parameter vector $\boldsymbol{\theta}_f$:
\begin{equation}
\min_{\boldsymbol{\theta}_f \in \mathbb{R}^{N_f}} \; \tilde{J}(\boldsymbol{\theta}_f) \coloneqq J\bigl(f(\boldsymbol{\theta}_f)\bigr) = \frac{1}{2} \, \bigl\| \mathcal{O}(\mathbf{u}(\boldsymbol{\theta}_f)) - \mathbf{y} \bigr\|_2^2 + \mathcal{R}(\boldsymbol{\theta}_f) \quad \text{s.t.} \quad \mathbf{R}^n = \mathbf{0}, \; \forall n.
\label{eq:discrete_optimization_problem}
\end{equation}

\paragraph{Discrete adjoint equation and reduced gradient} 
To derive the discrete adjoint equations, we introduce the Lagrangian function $\mathcal{L}$ by coupling the discrete objective function \eqref{eq:discrete_optimization_problem} with the sequence of forward residual constraints via the discrete adjoint trajectory $\boldsymbol{\lambda} = (\boldsymbol{\lambda}^1, \boldsymbol{\lambda}^2, \dots, \boldsymbol{\lambda}^{N_t})$, where $\boldsymbol{\lambda}^n \in \mathbb{R}^{N_x}$ denotes the Lagrange multiplier vector (or adjoint state) at the $n$-th time step:
\begin{equation}
\mathcal{L}\bigl(\mathbf{u}, \boldsymbol{\theta}_f, \boldsymbol{\lambda}\bigr) = \tilde{J}(\boldsymbol{\theta}_f) + \sum_{n=1}^{N_t} \bigl(\boldsymbol{\lambda}^{n}\bigr)^{\top} \mathbf{R}^n\bigl(\mathbf{u}^n, \mathbf{u}^{n-1}; \boldsymbol{\theta}_f\bigr).
\label{eq:lagrangian}
\end{equation}
By the principle of duality, the total derivative of the objective function with respect to the parameters matches the partial derivative of the Lagrangian ($\mathrm{d}_{\boldsymbol{\theta}_f} \tilde{J} = \partial_{\boldsymbol{\theta}_f} \mathcal{L}$) provided that the variations with respect to the state trajectory $\mathbf{u}$ vanish ($\partial_{\mathbf{u}} \mathcal{L} = \mathbf{0}$). Differentiating \eqref{eq:lagrangian} with respect to the state component $\mathbf{u}^{n}$ at each time step yields the adjoint state equation, which satisfies the backward-in-time linear recursion:

\begin{equation}
\Bigl(\partial_{\mathbf{u}^n} \mathbf{R}^n \Bigr)^{\!\top} \! \boldsymbol{\lambda}^n = -\Bigl(\partial_{\mathbf{u}^n} \mathbf{R}^{n+1} \Bigr)^{\!\top} \! \boldsymbol{\lambda}^{n+1} - \partial_{\mathbf{u}^n} \tilde{J}(\boldsymbol{\theta}_f), \qquad n=N_t-1, \dots, 1,
\label{eq:adjoint_recursion}
\end{equation}
where the terminal adjoint state $\boldsymbol{\lambda}^{N_t}$ is explicitly initialized at the final time step ($n=N_t$) by solving:
\begin{equation}
\Bigl(\partial_{\mathbf{u}^{N_t}} \mathbf{R}^{N_t} \Bigr)^{\!\top} \! \boldsymbol{\lambda}^{N_t} = - \partial_{\mathbf{u}^{N_t}} \tilde{J}(\boldsymbol{\theta}_f).
\label{eq:adjoint_terminal}
\end{equation}
The reduced gradient is then assembled as:
\begin{equation}
\nabla_{\boldsymbol{\theta}_f} \tilde{J}(\boldsymbol{\theta}_f) = \sum_{n=1}^{N_t} \Bigl(\partial_{\boldsymbol{\theta}_f} \mathbf{R}^n \Bigr)^{\!\top} \! \boldsymbol{\lambda}^n + \nabla_{\boldsymbol{\theta}_f} \mathcal{R}(\boldsymbol{\theta}_f).
\label{eq:reduced_gradient}
\end{equation}

We remark that for steady state problems, the adjoint equation \eqref{eq:adjoint_recursion} reduces to a single linear adjoint solve with the transpose of the full Jacobian of the converged forward state, and \eqref{eq:reduced_gradient} is computed by a single matrix--vector product.

\paragraph{Gradient verification}
The correctness of every discrete adjoint is established before any inversion:
for each benchmark, the reduced gradient $\nabla_{\boldsymbol{\theta}_f}
\tilde{J}$ is checked against centered finite differences of the discrete
objective along multiple randomly sampled perturbation directions, agreeing to
a maximum relative error below $10^{-3}$ across all tests. This verification guarantees that the measured differences between the two methods reflect algorithmic behavior rather than an erroneous sensitivity. \revaa{This verification is performed only once for each benchmark, prior to the inversion, as an implementation check; it is not part of the optimization loop and its cost is not included in the reported wall-clock timings.} 
 
\paragraph{Optimizer}
All adjoint inverse problems are solved with a second-order quasi-Newton method, SSBroyden.

 \subsection{Physics-Informed Neural Networks}
\label{subsec:pinn}

Physics-informed neural networks \cite{raissi2019physics} recast \eqref{eq:abstract_inverse} as the joint training of two neural networks: one approximating the state $u$  and one approximating the unknown function $f$, with the PDE constraint \eqref{eq:abstract_pde} enforced through a soft penalty on its residual.
This subsection describes the architectures, the composite loss, the collocation strategy, and the optimizer used in this work. 

\paragraph{Network setup}
Two multilayer perceptrons (MLPs) are used: a state network $N_{\boldsymbol{\theta}_u}$ approximating the solution $u$, and a parameter network $N_{\boldsymbol{\theta}_f}$ approximating the unknown function $f$, where $\boldsymbol{\theta}_u \in \mathbb{R}^{N_{u}}$ and $\boldsymbol{\theta}_f \in \mathbb{R}^{N_{f}}$ denote the vectors of trainable weights and biases for the respective networks.

\paragraph{Composite physics-informed loss}
The training objective is a weighted sum of mean-squared residuals enforcing the PDE, its boundary and initial conditions, and the observation data:
\begin{equation}
\ell(\boldsymbol{\theta}_u, \boldsymbol{\theta}_f) = w_{\mathrm{pde}} \ell_{\mathrm{pde}}(\boldsymbol{\theta}_u, \boldsymbol{\theta}_f) + w_{\mathrm{bc}} \ell_{\mathrm{bc}}(\boldsymbol{\theta}_u) + w_{\mathrm{ic}} \ell_{\mathrm{ic}}(\boldsymbol{\theta}_u) + w_{\mathrm{data}} \ell_{\mathrm{data}}(\boldsymbol{\theta}_u) + w_{\mathrm{reg}} \mathcal{R}(\boldsymbol{\theta}_f),
\label{eq:pinn_loss}
\end{equation}
with the first four terms defined by
\begin{align}
\ell_{\mathrm{pde}}(\boldsymbol{\theta}_u, \boldsymbol{\theta}_f) &= \frac{1}{N_{\mathrm{pde}}}\sum_{i=1}^{N_{\mathrm{pde}}} \Bigl|\partial_t N_{\boldsymbol{\theta}_u}+\mathcal{N}\bigl(N_{\boldsymbol{\theta}_u};N_{\boldsymbol{\theta}_f}\bigr)\Bigr|^2_{(\mathbf{x}_i,t_i)}, \label{eq:pinn_loss_pde}\\[2pt]
\ell_{\mathrm{bc}}(\boldsymbol{\theta}_u) &= \frac{1}{N_{\mathrm{bc}}}\sum_{i=1}^{N_{\mathrm{bc}}} \bigl|\mathcal{B}\bigl(N_{\boldsymbol{\theta}_u}\bigr)\bigr|^2_{(\mathbf{x}_i,t_i)}, \label{eq:pinn_loss_bc}\\[2pt]
\ell_{\mathrm{ic}}(\boldsymbol{\theta}_u) &= \frac{1}{N_{\mathrm{ic}}}\sum_{i=1}^{N_{\mathrm{ic}}} \bigl|N_{\boldsymbol{\theta}_u}(\mathbf{x}_i,0)-u_0(\mathbf{x}_i)\bigr|_2^2, \label{eq:pinn_loss_ic}\\[2pt]
\ell_{\mathrm{data}}(\boldsymbol{\theta}_u) &= \frac{1}{N_{\mathrm{data}}}\sum_{i=1}^{N_{\mathrm{data}}} \bigl|\mathcal{O}\bigl(N_{\boldsymbol{\theta}_u}\bigr)_i - y_i\bigr|_2^2. \label{eq:pinn_loss_data}
\end{align}
All spatial and temporal derivatives required to evaluate the residuals in \eqref{eq:pinn_loss_pde}--\eqref{eq:pinn_loss_bc} are obtained by automatic differentiation \cite{baydin2018automatic} of $N_{\boldsymbol{\theta}_u}$, so that $\ell$ is differentiable with respect to the network weights $(\boldsymbol{\theta}_u,\boldsymbol{\theta}_f)$.
The regularization term $\mathcal{R}(\boldsymbol{\theta}_f)$ matches the one used in the adjoint formulation \eqref{eq:abstract_inverse}, and all other loss weights are fixed to one in the experiments reported here so that no manual loss balancing is performed. \revaa{This choice is deliberate: equal weighting avoids artificially favoring any individual component of the physics-informed loss and eliminates benchmark-dependent loss-weight tuning from the comparison. Since all loss components are formulated as mean-squared residuals, unit weights provide a natural neutral baseline. In addition, the second-order SSBroyden optimizer used here was sufficiently robust to achieve convergence without introducing additional manual loss balancing.}
\revzz{SSBroyden updates have been reported to reach markedly lower
PINN training errors than first-order and limited-memory alternatives, both in
systematic optimizer comparisons~\cite{urban2025unveiling,KIYANI2025118308} and in
subsequent curvature-aware treatments~\cite{jnini2026curvature}; an open
implementation is available in~\cite{bioli2026selfscaled}.}

\paragraph{Training strategy}
Both the state network $N_{\boldsymbol{\theta}_u}$ and parameter network $N_{\boldsymbol{\theta}_f}$ utilize $\tanh$ activation functions with Xavier-uniform weight initialization \cite{glorot2010understanding} and are evaluated using double-precision arithmetic. 
Both networks are trained jointly using the second-order quasi-Newton optimizer SSBroyden, identical to the adjoint method. 
To prevent the networks from overfitting to a static point cloud, a budget of $N_{\mathrm{pde}}$ collocation points is drawn uniformly at random over the spatio-temporal domain $\Omega \times (0,T)$ and resampled after each outer quasi-Newton loop. 
The sampling effort is concentrated on this PDE residual term \eqref{eq:pinn_loss_pde} as it plays the dominant role in enforcing physical consistency. 
During each loop, the mini-batch is kept strictly unchanged to satisfy the deterministic objective requirement fundamental to BFGS-type line-search and curvature updates.

\subsection{Error metrics}

To systematically assess and compare the performance of these two para\-digms, we define discrete relative error metrics evaluated directly on the forward-solver mesh against the ground-truth arrays. 
Let $\mathbf{f}^{\star} \in \mathbb{R}^{N_x}$ denote the discrete ground-truth parameter field sampled on the grid, and let $\mathbf{u}^{\star} \in \mathbb{R}^{N_x \times N_t}$ represent the high-fidelity reference state trajectory obtained by running the numerical forward solver under this true parameter configuration.

The reconstruction accuracy of the parameter field is quantified by the relative discrete Euclidean norm $\varepsilon_f$:
\begin{equation}
\varepsilon_f = \frac{\bigl\| \mathbf{f}_{\mathrm{recovered}} - \mathbf{f}^{\star} \bigr\|_2}{\bigl\| \mathbf{f}^{\star} \bigr\|_2},
\label{eq:metric_epsilon_f}
\end{equation}
where $\mathbf{f}_{\mathrm{recovered}}$ is the identified parameter vector mapped onto the solver grid coordinates. 
For the adjoint method, $\mathbf{f}_{\mathrm{recovered}}^{\mathrm{adj}}$ 
is obtained from the discrete optimization routine, which may be interpolated onto the forward-solver mesh if their grid representations differ; for the PINN formulation, $\mathbf{f}_{\mathrm{recovered}}^{\mathrm{PINNs}}$ is obtained by evaluating the trained parameter network $N_{\boldsymbol{\theta}_f}$ directly at the mesh coordinates.

Similarly, the accuracy of the state trajectory recovery is measured by the relative discrete error $\varepsilon_u$:
\begin{equation}
\varepsilon_u = \frac{\bigl\| \mathbf{u}_{\mathrm{recovered}} - \mathbf{u}^{\star} \bigr\|_2}{\bigl\| \mathbf{u}^{\star} \bigr\|_2},
\label{eq:metric_epsilon_u}
\end{equation}
where $\mathbf{u}_{\mathrm{recovered}}$ represents discrete state vector to be compared against the reference $\mathbf{u}^\star$. 
For the adjoint method, $\mathbf{u}_{\mathrm{recovered}}^{\mathrm{adj}}$ is the discrete state vector computed during the optimization process. 
For the PINN approach, to ensure a \revzz{controlled} and consistent comparison, $\mathbf{u}_{\mathrm{recovered}}^{\mathrm{PINNs}}$ is the re-simulated state trajectory obtained by passing the parameter network's grid predictions $\mathbf{f}_{\mathrm{recovered}}^{\mathrm{PINNs}}$ into the discrete forward solver rather than evaluating the neural surrogate network $N_{\boldsymbol{\theta}_u}$. \revaa{This choice isolates the quality of the identified parameter from the approximation error of the auxiliary PINN state network: both recovered parameters are propagated through the same discrete forward model, so that $\varepsilon_u$
 provides a consistent cross-method measure of the physical state induced by the inversion result. The neural state reconstruction itself is considered separately where relevant, as in Test 4.}

\section{Benchmark Inverse Problems}
\label{sec:benchmarks}

This section instantiates the abstract inverse problem of
Section~\ref{sec:inverse} on four benchmark PDEs of increasing nonlinearity and dimensionality, listed in Table~\ref{tab:benchmark_overview}. 
We specify the forward operator $\mathcal{N}$, the boundary
operator $\mathcal{B}$, the initial condition $u_0$ and the unknown
function $f$ in the notation of~\eqref{eq:abstract_pde}, identify the
observation operator $\mathcal{O}$, instantiate the regularized inverse
objective~\eqref{eq:abstract_inverse}, and report the recovered fields.
Only the choices that affect the interpretation of the inverse-recovery quality are stated in the main text, whereas implementation details can be found in ~\ref{sec:appendix}. 
The general outcome of our test cases is summarized in Table~\ref{tab:inverse_summary}.
\revzz{Three conventions apply throughout. The discretization floor
$\varepsilon_u^{\mathrm{solver}}$ is the relative $L^2$ error of the forward solver
against a converged reference. Test~1 is quoted for
the neural-network parameterization and Test~2 at the $1\%$ noise level.}
All wall-clock timings reported in this paper were measured on a single compute node equipped with an NVIDIA L40S GPU and a dual-socket AMD EPYC 9554 CPU. No parallel computing is used.

\begin{table}[ht]
\centering
\footnotesize
\caption{Overview of the four inverse benchmarks.}
\label{tab:benchmark_overview}
\begin{tabular}{clcc}
\toprule
Test & Governing PDE & Unknown $f$ & Observation \\
\midrule
1 & unsteady 1D viscous Burgers     & spatial forcing $f(x)$                  & terminal state \\
2 & steady 2D Darcy                 & log-permeability $f(\mathbf{x})$        & sparse noisy probes \\
3 & unsteady 3D Allen--Cahn         & state-dependent forcing $f(u)$          & terminal state \\
4 & unsteady 2D Navier--Stokes      & scalar viscosity $\nu$                  & sparse wake probes \\
\bottomrule
\end{tabular}
\end{table}

\begin{table}[!htbp]
\centering
\footnotesize
\setlength{\tabcolsep}{3pt}
\caption{Summary of inverse-recovery quality and wall-clock cost across the four
benchmarks. \revzz{Entries are mean\,$\pm$\,standard deviation over five seeds;
$^{\dagger}$ marks deterministic estimators, reported from a single run.}}
\label{tab:inverse_summary}
\revzz{(a) Accuracy}\\[2pt]
\begin{tabular}{cccccc}
\toprule
Test
& $\varepsilon_f^{\mathrm{adj}}$
& $\varepsilon_f^{\mathrm{PINN}}$
& $\varepsilon_u^{\mathrm{solver}}$
& $\varepsilon_u^{\mathrm{adj}}$
& $\varepsilon_u^{\mathrm{PINN}}$ \\
\midrule
1 & \revzz{$3.09{\pm}3.4{\times}10^{-4}$}  & \revzz{$\mathbf{1.04{\pm}0.06{\times}10^{-4}}$}
  & $4.19{\times}10^{-5}$
  & \revzz{$\mathbf{1.64{\pm}1.5{\times}10^{-5}}$} & \revzz{$5.69{\pm}0.03{\times}10^{-5}$} \\
2 & \revzz{$2.87{\pm}0.13{\times}10^{-1}$} & \revzz{$\mathbf{2.70{\pm}0.18{\times}10^{-1}}$}
  & $7.45{\times}10^{-6}$
  & \revzz{$\mathbf{9.65{\pm}1.1{\times}10^{-3}}$} & \revzz{$1.00{\pm}0.06{\times}10^{-2}$} \\
3 & \revzz{$\mathbf{8.30{\pm}6.6{\times}10^{-4}}$} & \revzz{$3.88{\pm}2.2{\times}10^{-3}$}
  & $4.40{\times}10^{-4}$
  & \revzz{$\mathbf{9.13{\pm}9.5{\times}10^{-5}}$} & \revzz{$1.30{\pm}0.62{\times}10^{-3}$} \\
4 & \revzz{$\mathbf{6.83{\times}10^{-3}}^{\dagger}$} & \revzz{$6.80{\pm}4.8{\times}10^{-2}$}
  & \revzz{$1.41{\times}10^{-3}$}
  & \revzz{$\mathbf{1.65{\times}10^{-3}}^{\dagger}$} & \revzz{$6.63{\pm}4.2{\times}10^{-3}$} \\
\bottomrule
\end{tabular}

\vspace{6pt}
\revzz{(b) Cost: wall-clock, speedup and peak memory [GB]}\\[2pt]
\revzz{
\begin{tabular}{ccccccc}
\toprule
Test & $t_{\mathrm{adj}}$ [s] & $t_{\mathrm{PINN}}$ [s]
& $t_{\mathrm{adj}}/t_{\mathrm{PINN}}$
& adjoint & PINN, CPU & PINN, GPU \\
\midrule
1 & $1107{\pm}21$        & $\mathbf{195{\pm}4}$ & $5.7$   & $1.16$          & $1.23$          & $0.16$ \\
2 & $\mathbf{1.9{\pm}0.2}$ & $196{\pm}1$      & $0.010$ & $0.49$ & $1.55$          & $0.08$ \\
3 & $5825{\pm}350$       & $\mathbf{410{\pm}13}$ & $14.2$ & $5.59$          & $3.03$ & $1.28$ \\
4 & $1722^{\dagger}$     & $\mathbf{196{\pm}6}$ & $8.8$   & $1.78$          & $1.22$ & $0.87$ \\
\bottomrule
\end{tabular}
}
\end{table}

\subsection{Test 1: Linear forcing identification for the 1D viscous Burgers equation}
\label{subsec:burgers1d_inverse}

The first benchmark instantiates~\eqref{eq:abstract_pde} as the
one-dimensional viscous Burgers equation on $\Omega=[0,L)$ with $L=2\pi$ and
time interval $[0,T]$ with $T=1$. Periodic boundary conditions are used so
that the inverse problem is driven entirely by the forcing $f(x)$, governed by the system:
\begin{equation}
\begin{aligned}
&\partial_t u + u\,\partial_x u - \nu\,\partial_{xx}u = f(x), & x\in[0,L),\ t\in(0,T],\\
&u(0,t)=u(L,t), \partial_x u(0,t)=\partial_x u(L,t), & t\in(0,T], \\
&u(x,0)=\sin(x), & x\in[0,L).
\end{aligned}
\label{eq:burgers1d}
\end{equation}
The unknown function is the stationary forcing $f(x)$, with reference profile
$f^\star(x) = \sin(2x)$. The forward problem is discretized in space by
second-order centered finite differences on a uniform periodic grid of
$N_x = 512$ nodes and in time by a Crank--Nicolson scheme.
The observation $\mathbf{y} \in \mathbb{R}^{N_x}$ is the terminal
state $\mathbf{u}^{N_t}$ produced by the same solver under the reference
forcing $f^\star$ sampled on the forward grid ($N_x = 512$).
The objective function in the inverse problem~\eqref{eq:abstract_inverse} reduces to
\begin{equation}
\tilde{J}(\boldsymbol{\theta}_f)
\;=\;\frac{\Delta x}{2}\,
\bigl\|\mathbf{u}^{N_t}(\boldsymbol{\theta}_f)-\mathbf{y}\bigr\|_2^{2},
\label{eq:burgers1d_objective}
\end{equation}
with no explicit regularization ($\mathcal{R}\equiv 0$).

\revzz{To separate the effect of the representation from that of the inversion
algorithm, this benchmark is run over two methods and two
representations, with everything else held fixed.} The two
parameterizations are the neural field
$f(\boldsymbol{\theta}_f)=N_{\boldsymbol{\theta}_f}(x)$, an MLP with two hidden
layers of width $32$ and $\tanh$ activations ($N_f=1153$), and a coarse mesh in
which $f$ is represented by $N_f=64$ degrees of freedom on a uniform
piecewise-linear mesh and linearly interpolated to the forward grid. \revzz{Both
parameterizations are shared exactly: the two methods use the same MLP, and
the PINN evaluates the same piecewise-linear prolongation as the adjoint, pointwise
at its collocation points.
Each cell is averaged over five random initializations, except adjoint\,+\,grid,
which is deterministic.}

\begin{figure}[!htbp]
    \centering
    \includegraphics[width=\textwidth]{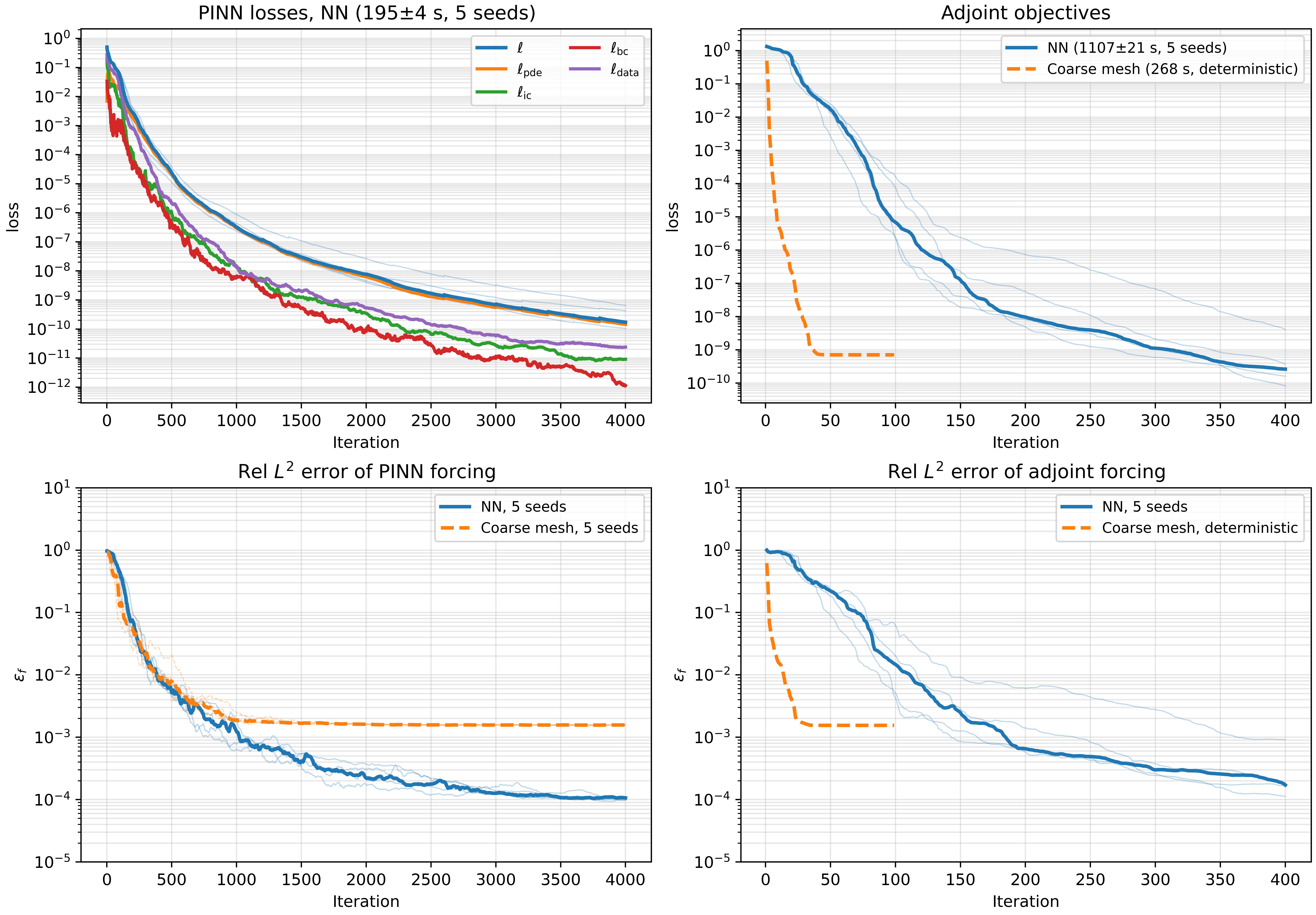}
    \caption{{\bf Test 1 - 1D unsteady Burgers}. Convergence histories. Horizontal
    axes are quasi-Newton iterations. \revzz{Thin curves are five initializations and bold is their median. The coarse-grid adjoint is deterministic and single.}}
    \label{fig:burgers_history_comparison}
\end{figure}

\begin{figure}[!htbp]
    \centering
    \begin{subfigure}[b]{0.85\textwidth}
        \centering
        \includegraphics[width=\textwidth]{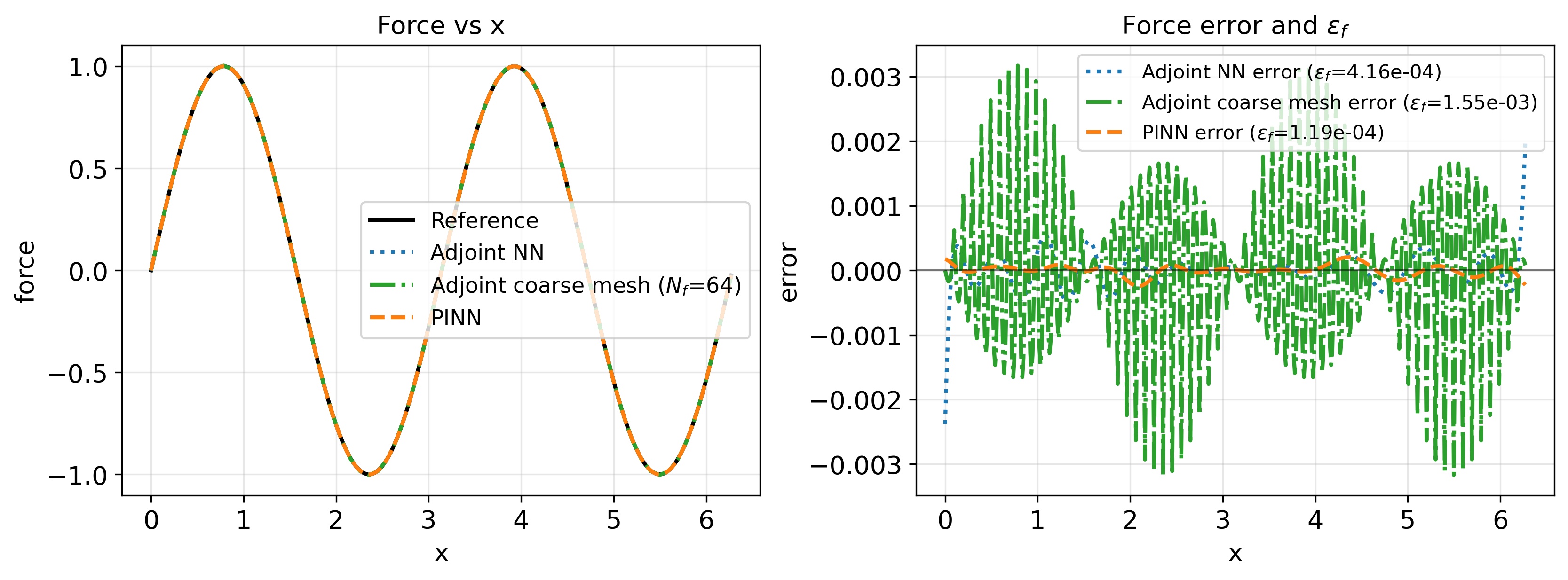}
        \caption{Recovered forcing $\mathbf{f}_{\mathrm{recovered}}$
        and pointwise error $\mathbf{f}_{\mathrm{recovered}}-\mathbf{f}^\star$.}
        \label{fig:burgers_force_comparison}
    \end{subfigure}\\[0.5em]
    \begin{subfigure}[b]{0.85\textwidth}
        \centering
        \includegraphics[width=\textwidth]{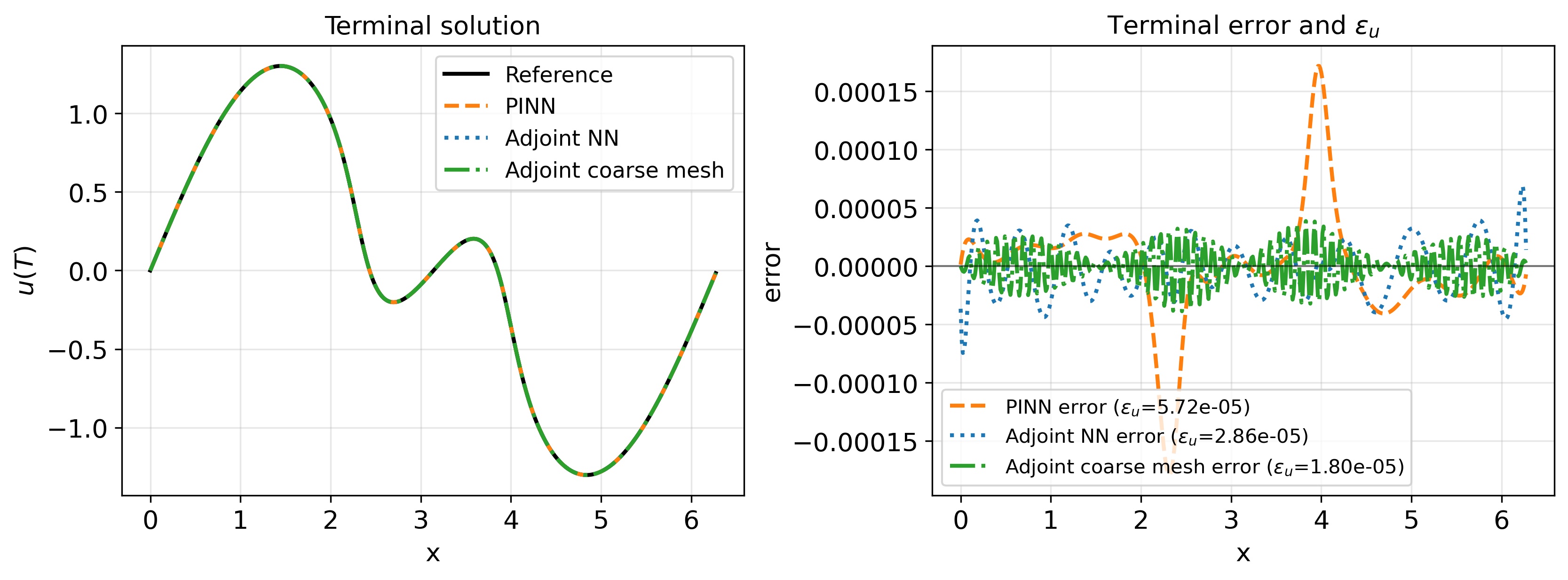}
        \caption{Recovered terminal state $\mathbf{u}_{\mathrm{recovered}}^{N_t}$
        and pointwise error $\mathbf{u}^{N_t}_{\mathrm{recovered}}-{\mathbf{u}^\star}^{N_t}$.}
        \label{fig:burgers_solution_comparison}
    \end{subfigure}
    \caption{{\bf Test 1 - 1D unsteady Burgers}. Recovered forcing field (top) and terminal-state field (bottom). $\varepsilon_u$ is calculated based on the terminal state $\mathbf{u}^{N_t}$. \revzz{Single initialization.}}
    \label{fig:burgers_force_solution}
\end{figure}

The convergence histories are reported in
Fig.~\ref{fig:burgers_history_comparison} and the recovered fields in
Fig.~\ref{fig:burgers_force_solution}.
All estimators recover the spatial forcing accurately with $\varepsilon_f<2\times 10^{-3}$ and reproduce the terminal solution to within $10^{-4}$.

\revzz{The four combinations, summarized in Table~\ref{tab:repr_factorial},
separate the two effects. Moving from the coarse grid to the neural field improves the parameter
error by roughly an order of magnitude for both algorithms, from
$1.5\text{--}1.6\times10^{-3}$ to $1.0\text{--}3.1\times10^{-4}$, whereas
exchanging the algorithm at a fixed representation changes it far less: on the grid
the two methods are statistically indistinguishable. The representation is
therefore the dominant factor and the algorithm the secondary one.}

\revzz{The two representations therefore suit the two methods differently. The
coarse grid is the natural setting for the adjoint: it cuts its cost fourfold,
from $1107\pm21$\,s to $268$\,s, and drives the state error to
$1.80\times10^{-5}$, statistically tied with the best of the four and at the
discretization floor, because the discrete adjoint extracts as much state-matching
accuracy as the forward discretization allows. What it cannot do is
resolve the unknown beyond its basis --- both methods stall at
$\varepsilon_f\approx1.6\times10^{-3}$ on that mesh, a ceiling belonging to the
piecewise-linear basis rather than to either algorithm. On the neural field the
ordering reverses: the PINN is both more accurate on the forcing
($1.04\pm0.06\times10^{-4}$ against $3.09\pm3.4\times10^{-4}$) and roughly five
times faster ($195\pm4$\,s against $1107\pm21$\,s), because it never integrates
the discrete trajectory, whereas each adjoint iteration sweeps forward and
backward through it and optimizes in a $10^3$-dimensional weight space. The
adjoint in weight space is also the less repeatable of the two: its five seeds of $\varepsilon_f$
cluster near $1.7\times10^{-4}$ but one reaches $9.1\times10^{-4}$, a largest-to-smallest ratio of $8.1$ against the PINN's $1.1$. Rather than degrading continuously, it occasionally gets trapped in a significantly worse local minimum.}

\begin{table}[ht]
\centering
\color{zzrev}
\footnotesize
\setlength{\tabcolsep}{5pt}
\caption{\revzz{{\bf Test 1 --- two methods $\times$ two representations.} Mean\,$\pm$\,standard deviation
over five initializations; $^{\dagger}$ deterministic, single run. For both methods
$\varepsilon_u$ is obtained by pushing the recovered forcing back through the
discrete solver.}}
\label{tab:repr_factorial}
\begin{tabular}{llccc}
\toprule
Representation & Algorithm & $\varepsilon_f$ & $\varepsilon_u$ & $t$ [s] \\
\midrule
Coarse grid, $N_f=64$   & Adjoint & $1.549{\times}10^{-3}$$^{\dagger}$          & $1.805{\times}10^{-5}$$^{\dagger}$ & $268$$^{\dagger}$ \\
Coarse grid, $N_f=64$   & PINN    & $1.574{\pm}0.010{\times}10^{-3}$            & $1.007{\pm}0.098{\times}10^{-4}$ & $\mathbf{69{\pm}5}$ \\
\midrule
Neural field, $N_f=1153$ & Adjoint & $3.09{\pm}3.4{\times}10^{-4}$              & $\mathbf{1.64{\pm}1.5{\times}10^{-5}}$ & $1107{\pm}21$ \\
Neural field, $N_f=1153$ & PINN    & $\mathbf{1.04{\pm}0.06{\times}10^{-4}}$     & $5.69{\pm}0.03{\times}10^{-5}$   & $195{\pm}4$ \\
\bottomrule
\end{tabular}
\end{table}

\subsection{Test 2: Linear permeability identification for 2D Darcy flow under noisy observations}
\label{subsec:darcy_inverse}

The second benchmark instantiates~\eqref{eq:abstract_pde} as the steady
Darcy equation on the unit square $\Omega=(0,1)^2$ with homogeneous
Dirichlet boundary conditions and known constant volumetric source $g=100$, governed by the system:
\begin{equation}
\begin{aligned}
&-\nabla\!\cdot\!\bigl(k(\mathbf{x})\,\nabla u(\mathbf{x})\bigr) = g, & \mathbf{x}\in\Omega,\\
&u(\mathbf{x})=0, & \mathbf{x}\in\partial\Omega.
\end{aligned}
\label{eq:darcy}
\end{equation}
The unknown function is the log-permeability $f(\mathbf{x})$, which parameterizes the physical permeability as
$k(\mathbf{x})=\exp\bigl(f(\mathbf{x})\bigr)$ so that positivity of $k$ is enforced by construction.
The reference $f^\star$ is a sample from a Gaussian random field with $r=128$ Karhunen--Lo\`eve modes
(\ref{subsec:darcy_implementation}). The forward problem is discretized by continuous
piecewise quadratic $\mathbb{P}_2$ finite elements on a structured triangular mesh of $32\times 32$ cells,
each split into two triangles, giving $N_x=4225$ degrees of freedom.

\begin{figure}[!htbp]
    \centering
    \includegraphics[width=\textwidth]{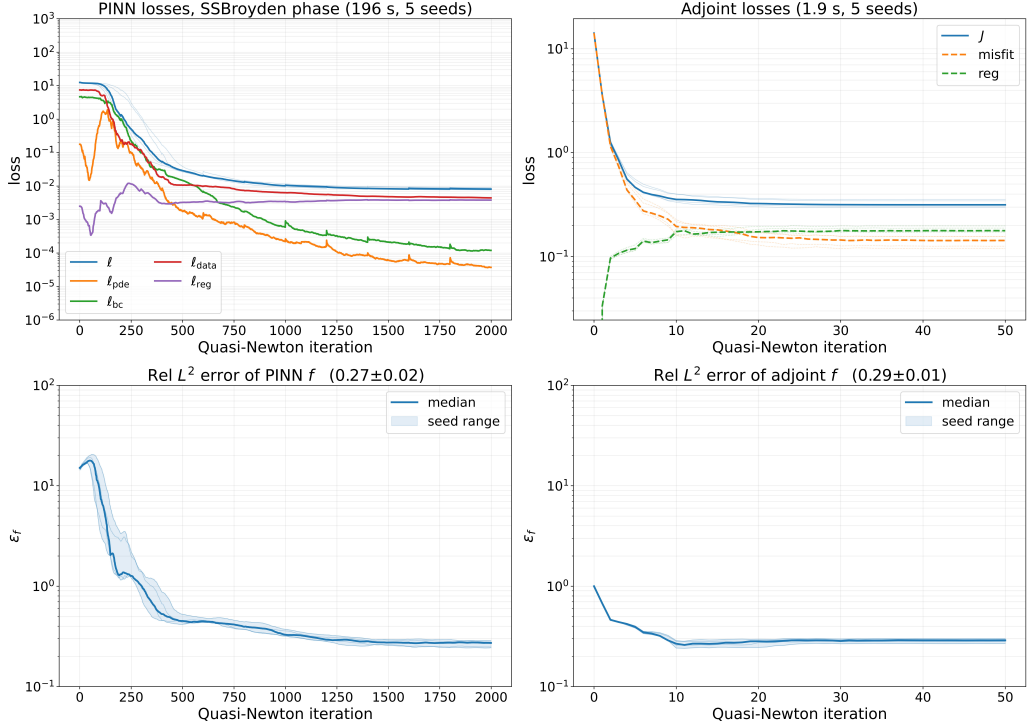}
    \caption{{\bf Test 2 - 2D Darcy.} Convergence histories with noisy observations;
    horizontal axes are quasi-Newton iterations. \revzz{At $\sigma=1\%$ and
    $\gamma^{\ast}=10^{-2}$, the setting of Tables~\ref{tab:inverse_summary}
    and~\ref{tab:darcy_noise}. Thin curves are five seeds, bold their median. The
    seed varies the noise draw as well as the initialization, so the adjoint band,
    which is otherwise deterministic, is entirely data noise.}}
    \label{fig:darcy_training_history}
\end{figure}

\begin{figure}[!htbp]
    \centering
    \includegraphics[width=\textwidth]{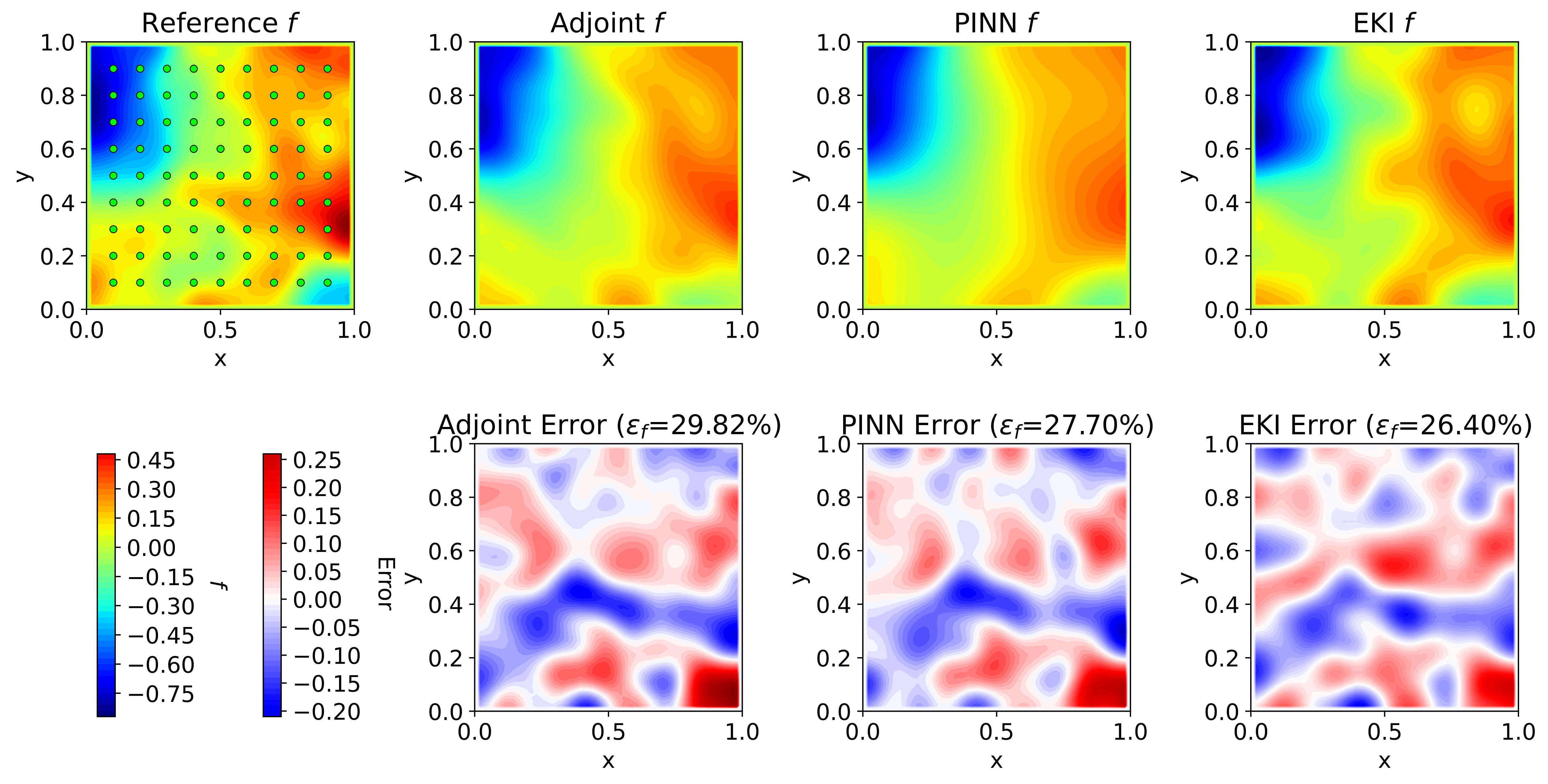}
    \caption{{\bf Test 2 - 2D Darcy.} Recovered log-permeability $\mathbf{f}_{\mathrm{recovered}}$ (top) and pointwise error $\mathbf{f}_{\mathrm{recovered}}-\mathbf{f}^\star$ (bottom), shown for the adjoint, PINN, and EKI estimators. The $9\times 9$ observation grid is overlaid on the reference panel (green dots). \revzz{At $\sigma=1\%$ noise and $\gamma^{\ast}=10^{-2}$, the setting of Table~\ref{tab:inverse_summary}; single initialization.}}
    \label{fig:darcy_m_comparison}
\end{figure}

\begin{figure}[!htbp]
    \centering
    \includegraphics[width=\textwidth]{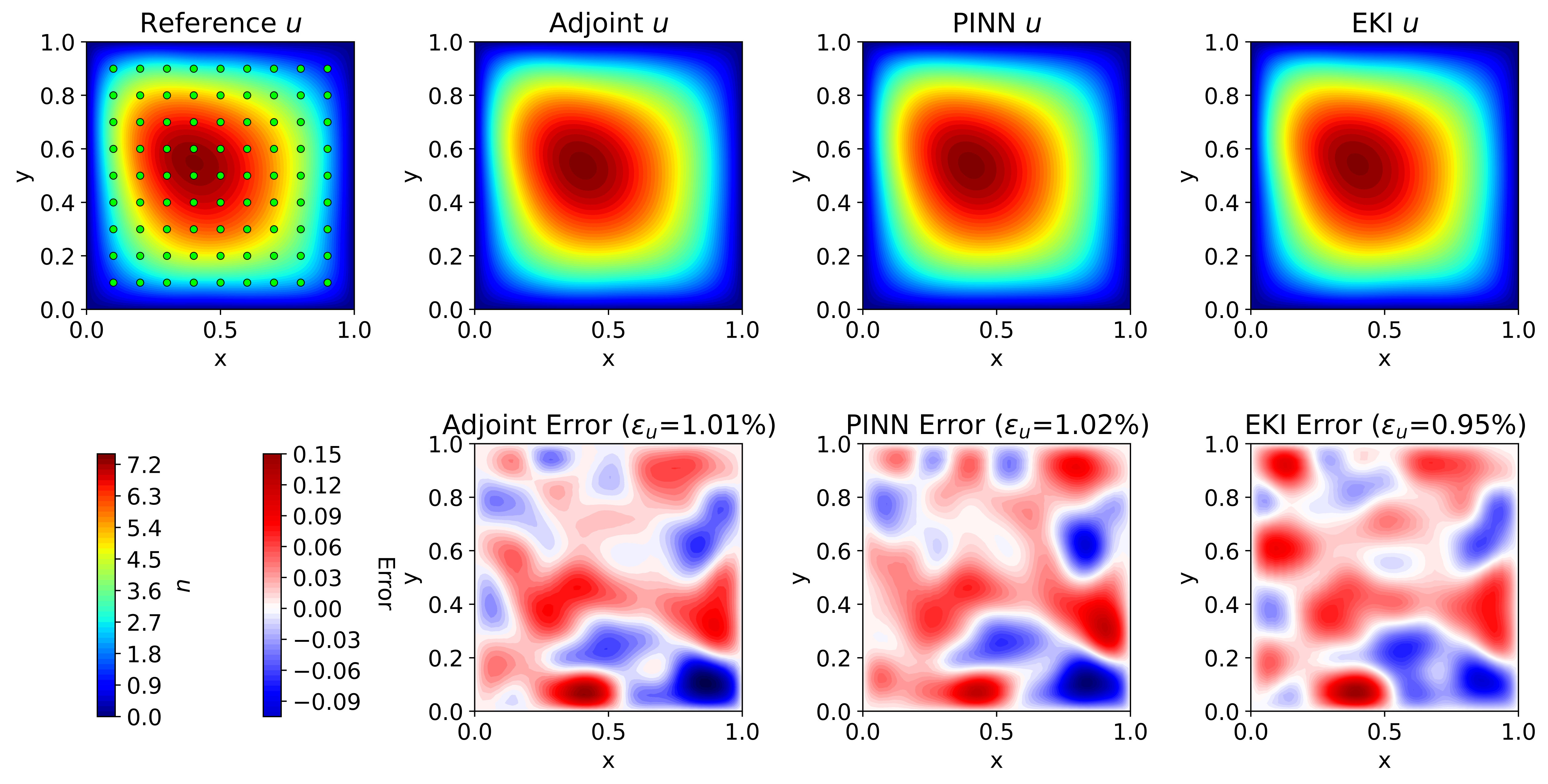}
    \caption{{\bf Test 2 - 2D Darcy.} Recovered state $\mathbf{u}_{\mathrm{recovered}}$ (top) and pointwise error $\mathbf{u}_{\mathrm{recovered}}-\mathbf{u}^\star$ (bottom). \revzz{At $\sigma=1\%$ noise and $\gamma^{\ast}=10^{-2}$; single initialization, so the errors on the panels are that one draw rather than the five-seed means quoted in the text and in Table~\ref{tab:inverse_summary}.}}
    \label{fig:darcy_u_comparison}
\end{figure}

The observation is sparse and noisy:
the data $\mathbf{y}\in\mathbb{R}^{81}$ are pointwise pressure measurements at a $9\times 9$ uniform
grid of interior locations, obtained by applying the exact $\mathbb{P}_2$ interpolation matrix
$P\in\mathbb{R}^{81\times N_x}$ to the state $\mathbf{u}^\star$ produced by the same solver
under the reference log-permeability $f^\star$ sampled on the forward mesh, and corrupted by $1\%$ additive Gaussian noise:
\begin{equation}
\mathbf{y}=P\,\mathbf{u}^\star+\boldsymbol{\eta},
\qquad
\boldsymbol{\eta}\sim\mathcal{N}\bigl(\mathbf{0},\sigma_\varepsilon^{\,2}\,I\bigr),
\qquad
\sigma_\varepsilon=10^{-2}\cdot\max_{i}|u^\star_i|.
\label{eq:darcy_obs}
\end{equation}
The regularized inverse objective~\eqref{eq:abstract_inverse} reduces to
\begin{equation}
\tilde{J}(\boldsymbol{\theta}_f)
\;=\;\frac{1}{2}\,
\bigl\|P\,\mathbf{u}(\boldsymbol{\theta}_f)-\mathbf{y}\bigr\|_2^{2}
\;+\;\mathcal{R}(\boldsymbol{\theta}_f),
\label{eq:darcy_objective}
\end{equation}
where $\mathcal{R}(f) = \tfrac{\gamma}{2}\int_{\Omega} \lvert \nabla f \rvert^2\,\mathrm{d}\mathbf{x}$
is an $H^1$ seminorm smoothness penalty of weight \revzz{$\gamma=\gamma^{\ast}(\sigma)$,
selected once per noise level and shared by both gradient-based methods as described
in Table~\ref{tab:darcy_noise}; at the $\sigma=1\%$ of the production runs
$\gamma^{\ast}=10^{-2}$},
discretized as detailed in~\ref{subsec:darcy_implementation}.

Three reconstructions are compared: the PINN with the state network
$N_{\boldsymbol{\theta}_u}(\mathbf{x})$ and the parameter network
$N_{\boldsymbol{\theta}_f}(\mathbf{x})$, both MLPs with three hidden layers of width
$32$ and $\tanh$ activations; the adjoint method optimizing directly in the
$128$-dimensional KL-coefficient space of the prior (``Adjoint KL''); and ensemble
Kalman inversion (``EKI'')~\cite{iglesias2013ensemble}, a derivative-free baseline run in the same
$128$-dimensional KL space. The adjoint and EKI estimators therefore share their
prior representation exactly and differ only in their inversion machinery.

The convergence histories of the two gradient-based estimators are reported in Fig.~\ref{fig:darcy_training_history}, whereas 
the recovered log-permeabilities in Fig.~\ref{fig:darcy_m_comparison}, and the corresponding states in
Fig.~\ref{fig:darcy_u_comparison}.
All three estimators recover the dominant large-scale structure of $f^\star$ and reproduce the steady state to within
$1\%$ relative error, but they differ in fine-scale details that the sparse noisy data cannot constrain.

The three log-permeability errors fall within a narrow band:
\revzz{$\varepsilon_f^{\mathrm{EKI}}=25.8\pm2.5\%$,
$\varepsilon_f^{\mathrm{PINN}}=27.0\pm1.8\%$ and
$\varepsilon_f^{\mathrm{Adj}}=28.7\pm1.3\%$, a spread of under three percentage
points that the seed-to-seed standard deviations do not resolve. At the state level
the same ordering holds, $\varepsilon_u^{\mathrm{EKI}}=0.90\pm0.13\%$,
$\varepsilon_u^{\mathrm{Adj}}=0.97\pm0.11\%$ and
$\varepsilon_u^{\mathrm{PINN}}=1.00\pm0.06\%$}, all three orders of
magnitude above the FEM discretization floor
$\varepsilon_u^{\mathrm{solver}}=7.45\times 10^{-6}$.

The qualitative differences visible in Fig.~\ref{fig:darcy_m_comparison} reflect the implicit prior of each method:
the PINN reconstruction is visibly smoother than both the Adjoint KL and EKI fields, consistent with the implicit
smoothness prior of the neural-field representation; the Adjoint KL and EKI fields exhibit sharper
local features inherited from the explicit KL truncation that both methods share.

The wall-clock costs span nearly two orders of magnitude. The deterministic
adjoint converges in \revzz{$1.9$\,s} by operating in a $128$-dimensional KL space on a
$\mathbb{P}_2$ FEM stiffness matrix that admits sparse direct factorization.
\revzz{EKI costs $164\pm20$\,s: each of its $5$ iterations evaluates the forward
solver once per ensemble member over $2000$ members. The PINN is the most
expensive at $196\pm1$\,s} because both $N_{\boldsymbol{\theta}_u}$ and
$N_{\boldsymbol{\theta}_f}$ are jointly represented and optimized as neural
networks over a $\sim\!10^{3}$-dimensional parameter space. 

\revzz{
\begin{table}[!htbp]
\centering
\color{zzrev}
\footnotesize
\caption{\revzz{{\bf Test 2 --- noise sensitivity.} $\varepsilon_f$ at the selected
$\gamma^{\ast}$, mean\,$\pm$\,standard deviation over five seeds; $^{\dagger}$ marks the
one deterministic entry, the adjoint at $\sigma=0$, where the seed varies only the
noise draw and there is none. $\gamma^{\ast}$ is
chosen once per noise level and applied to both gradient-based methods, since the
regularization weight belongs to the inverse problem rather than to the solver;
EKI carries no explicit $\gamma$, its regularization coming from the ensemble
prior. $\gamma^{\ast}$ is selected from the grid
$\{0,10^{-4},10^{-3},10^{-2},3{\times}10^{-2},10^{-1}\}$.}}
\label{tab:darcy_noise}
\begin{tabular}{ccccc}
\toprule
$\sigma$ & $\gamma^{\ast}$ & adjoint & PINN & EKI \\
\midrule
$0$     & $0$       & $0.087^{\dagger}$ & $0.191{\pm}0.008$ & $0.078{\pm}0.005$ \\
$0.1\%$ & $10^{-4}$ & $0.109{\pm}0.005$ & $0.186{\pm}0.007$ & $0.112{\pm}0.006$ \\
$1\%$   & $10^{-2}$ & $0.287{\pm}0.013$ & $0.270{\pm}0.018$ & $0.258{\pm}0.025$ \\
$5\%$   & $10^{-1}$ & $0.442{\pm}0.048$ & $0.444{\pm}0.043$ & $0.430{\pm}0.073$ \\
$10\%$  & $10^{-1}$ & $0.548{\pm}0.103$ & $0.547{\pm}0.082$ & $0.526{\pm}0.105$ \\
\bottomrule
\end{tabular}
\end{table}
}

\revzz{Table~\ref{tab:darcy_noise} sweeps the noise level and the regularization
weight, five seeds per cell, with $\gamma^{\ast}$ selected once per noise level and
shared by both gradient-based methods since it belongs to the inverse problem
rather than to the solver. From $1\%$ noise upwards the three estimators are
statistically indistinguishable: at $\sigma=5\%$ they give $0.442$, $0.444$ and
$0.430$ with standard deviations of $0.04$--$0.07$. What limits the recovered
permeability here is therefore the problem, not the method. The forward map
$f\mapsto u$ has a strongly contracting spectrum, so moderate differences in $f$
project onto comparable observation residuals and the sparse noisy probes cannot
distinguish among the reconstructions; once the noise sets the achievable
accuracy, all three estimators reach it, and the choice among them reduces to
cost.}

\subsection{Test 3: Nonlinear state-dependent forcing identification for the 3D Allen--Cahn equation}
\label{subsec:ac3d_inverse}

The third benchmark instantiates~\eqref{eq:abstract_pde} as the
three-dimensional Allen--Cahn equation on $\Omega=(0,L)^3$ with $L=1$,
interface width $\varepsilon=10^{-2}$, and time interval $[0,T]$ with $T=3$.
Homogeneous Neumann boundary
conditions are imposed on all faces, governed by the system:
\begin{equation}
\begin{aligned}
&\partial_t u - \varepsilon^{2}\bigl(\partial_{xx}u+\partial_{yy}u+\partial_{zz}u\bigr) = f(u), & \mathbf{x}\in\Omega,\ t\in(0,T],\\
&\nabla u\cdot\mathbf{n}=0, & \mathbf{x}\in\partial\Omega,\ t\in(0,T], \\
&u(\mathbf{x},0)=\cos(\pi x)\cos(\pi y)\cos(\pi z), & \mathbf{x}\in\Omega.
\end{aligned}
\label{eq:ac3d}
\end{equation}
The unknown function is a
\emph{state-dependent} reaction $f(u)$, a scalar nonlinear functional of the solution, with reference
profile $f^\star(u)=u-u^{3}$ generating bistable phase-field dynamics with a sharp transition layer.
This qualitatively changes the inverse problem: the forward map $f\mapsto u(f)$ is itself nonlinear in $f$,
and the long-time dissipative dynamics contract any component of the sensitivity orthogonal to the leading
attracting mode, so the conditioning of the terminal-state sensitivity deteriorates as $T$ grows.

The tensor-product cosine initial condition develops eight diffuse interfaces meeting at
the center of the cube, so the terminal state contains both well-separated bulk regions ($u\approx\pm 1$) and
sharp transition layers.

The forward problem is discretized in space by the second-order seven-point
finite-difference Laplacian with ghost-point Neumann closure on a uniform
tensor-product grid of $N_x=N_y=N_z=129$ nodes, and in time by a second-order
IMEX BDF2/AB2 scheme with $\Delta t=5\times 10^{-2}$. Because the implicit
operator is the constant-coefficient Neumann Laplacian, it is diagonalized
exactly by the discrete cosine basis~\cite{strang1999dct}, so each step is a
direct, matrix-free $\mathcal{O}(N\log N)$ transform--divide--inverse-transform
operation, incurring none of the fill-in of a three-dimensional sparse
factorization and reused unchanged for the adjoint transpose solve
(\ref{subsec:ac3d_implementation}). The observation $\mathbf{y}$ is
the terminal state $\mathbf{u}^{N_t}$ produced by the same solver under the
reference reaction $f^\star$ on the forward grid, and the inverse
problem~\eqref{eq:abstract_inverse} reduces to
\begin{equation}
\tilde{J}(\boldsymbol{\theta}_f)
\;=\;\frac{\Delta_x \Delta_y \Delta_z }{2}\,
\bigl\|\mathbf{u}^{N_t}(\boldsymbol{\theta}_f)-\mathbf{y}\bigr\|_2^{2},
\label{eq:ac3d_objective}
\end{equation}
with no explicit regularization ($\mathcal{R}\equiv 0$).

Two reconstructions are compared on this benchmark: the PINN with the state network
$N_{\boldsymbol{\theta}_u}(x,y,z,t)$ and the parameter network $N_{\boldsymbol{\theta}_f}(u)$;
and the adjoint method with the same MLP parameterization
$f(\boldsymbol{\theta}_f)=N_{\boldsymbol{\theta}_f}(u)$ (``Adjoint NN''). 
The parameter network is an MLP with two hidden layers of width
$32$ and $\tanh$ activations, while the state network has three hidden layers of
width $32$ and the same activation.

\begin{figure}[!htbp]
    \centering
    \includegraphics[width=\textwidth]{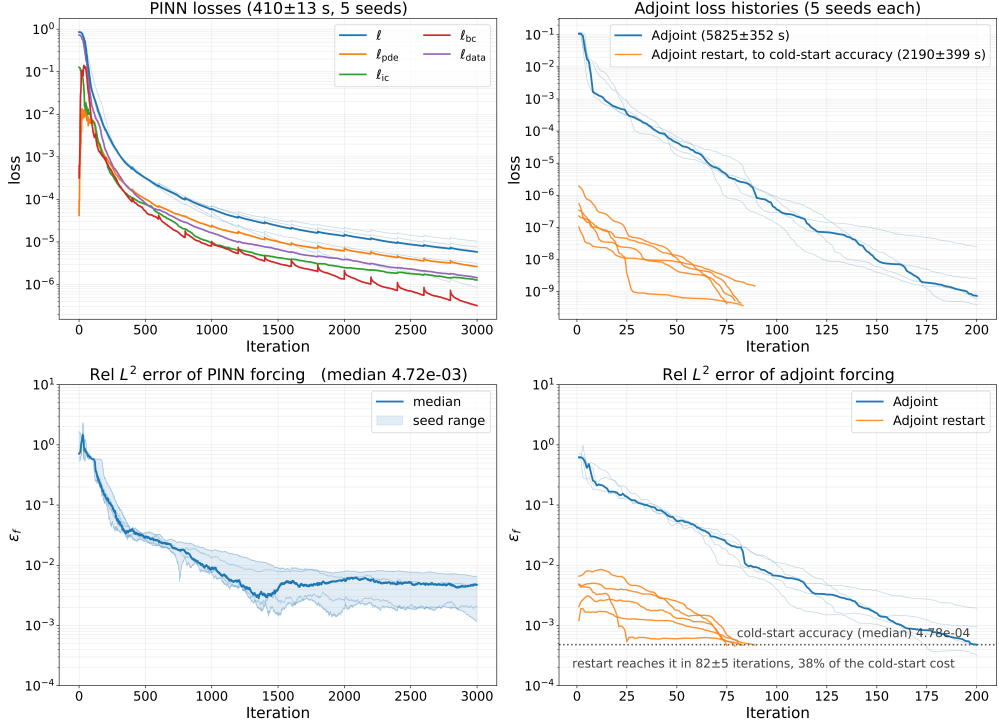}
    \caption{{\bf Test 3 - 3D Allen--Cahn.} Convergence histories; horizontal axes are
    accepted quasi-Newton iterations. \revzz{Thin curves are five seeds per arm,
    bold and the dotted target their medians. Each restart is drawn only as far as
    the reaction error the cold start reaches, after $82\pm5$ iterations against
    $200$, or $38\%$ of its wall-clock.}}
    \label{fig:ac3d_history_comparison}
\end{figure}

\begin{figure}[!htbp]
    \centering
    \includegraphics[width=\textwidth]{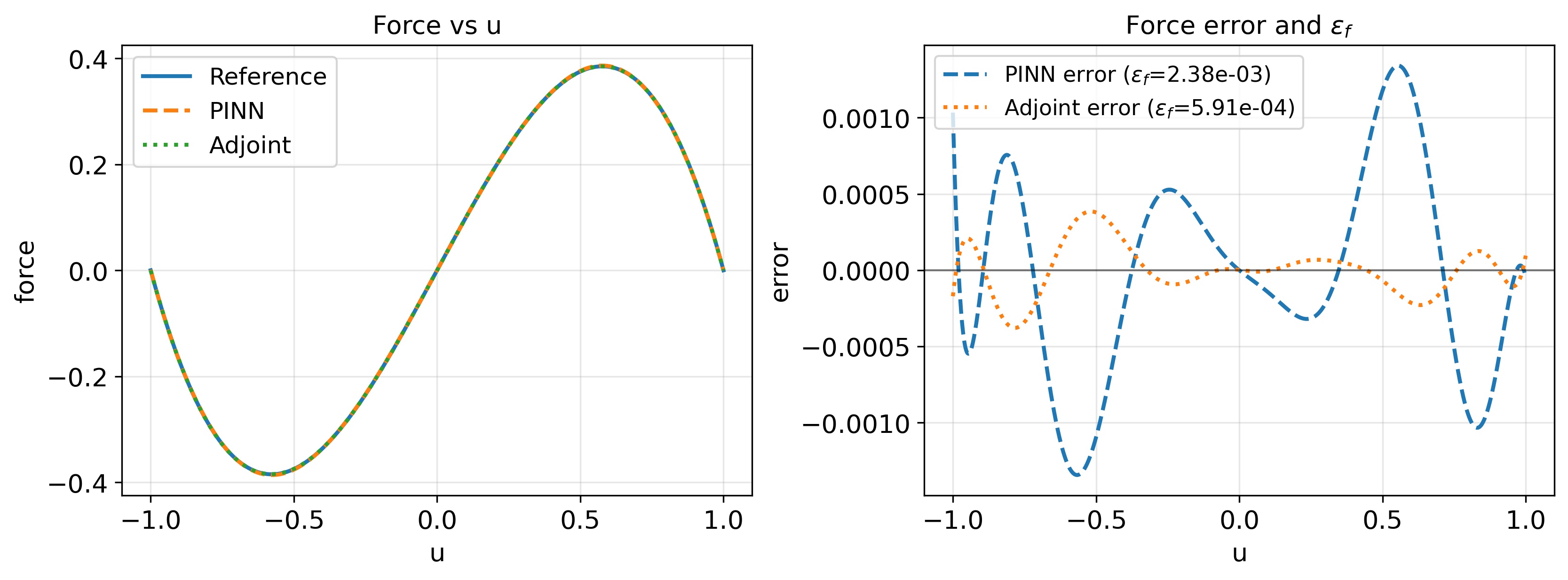}
    \caption{{\bf Test 3 - 3D Allen--Cahn.} Recovered reaction term $\mathbf{f}_{\mathrm{recovered}}$ (left) and pointwise error $\mathbf{f}_{\mathrm{recovered}}-\mathbf{f}^\star$ (right) on the physical range $u\in[-1,1]$. \revzz{Single initialization.}}
    \label{fig:ac3d_force_comparison}
\end{figure}

\begin{figure}[!htbp]
    \centering
    \includegraphics[width=\textwidth]{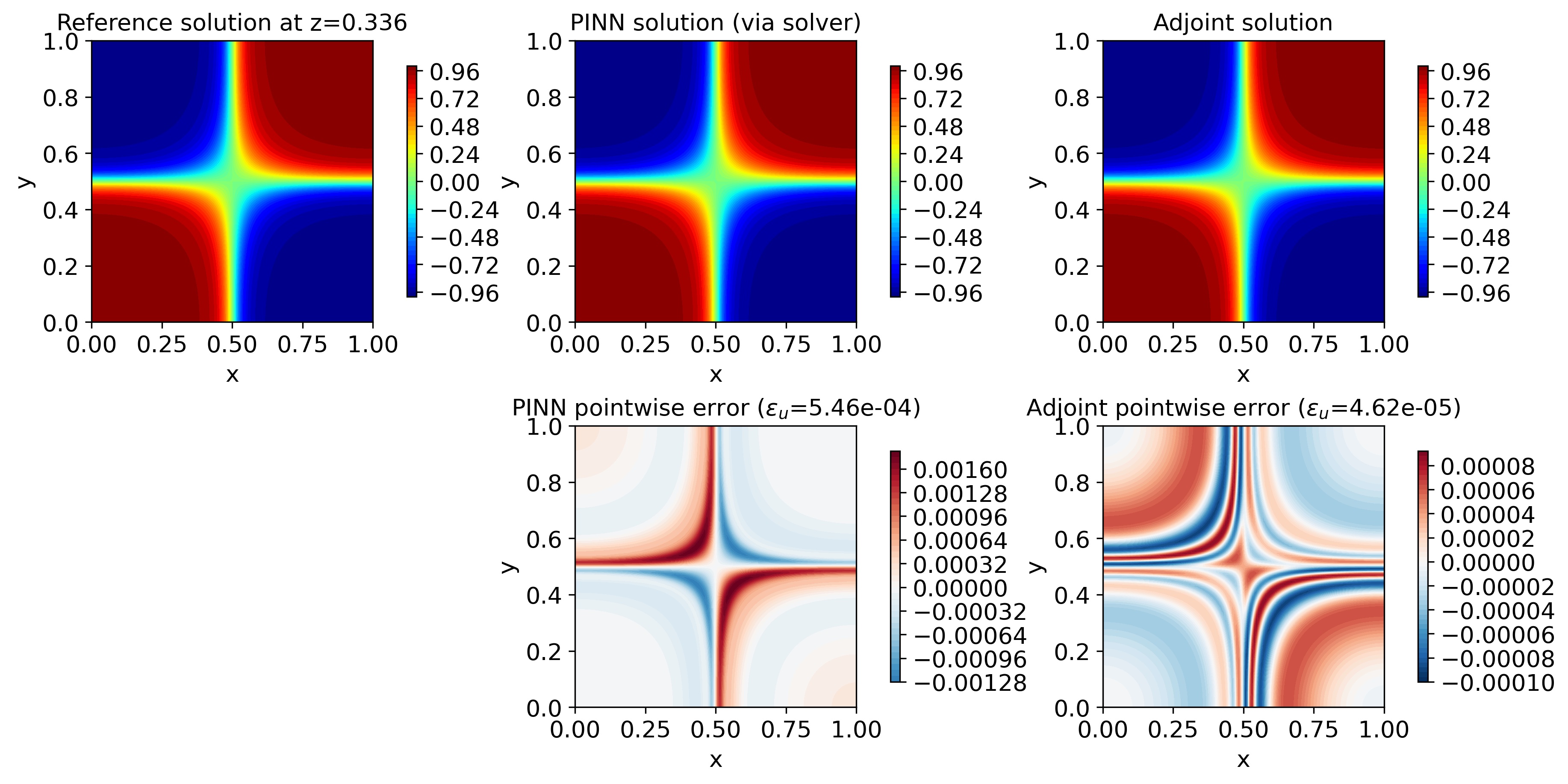}
    \caption{{\bf Test 3 - 3D Allen--Cahn.} Slice ($z=0.336$) of the recovered terminal state $\mathbf{u}_{\mathrm{recovered}}^{N_t}$ (top) and pointwise error relative to ${\mathbf{u}^\star}^{N_t}$ (bottom). $\varepsilon_u$ is calculated based on the full terminal state $\mathbf{u}^{N_t}$. \revzz{Single initialization.}}
    \label{fig:ac3d_solution_comparison}
\end{figure}

The convergence histories of both estimators are reported in
Fig.~\ref{fig:ac3d_history_comparison}, the recovered reaction in
Fig.~\ref{fig:ac3d_force_comparison}, and the recovered terminal state in
Fig.~\ref{fig:ac3d_solution_comparison}. The adjoint method is the more
accurate of the two on both metrics: \revzz{the parameter error is
$\varepsilon_f^{\mathrm{Adj}}=8.30\pm6.6\times 10^{-4}$ versus
$\varepsilon_f^{\mathrm{PINN}}=3.88\pm2.2\times 10^{-3}$, and the terminal-state
error is $\varepsilon_u^{\mathrm{Adj}}=9.13\pm9.5\times 10^{-5}$ versus
$\varepsilon_u^{\mathrm{PINN}}=1.30\pm0.62\times 10^{-3}$, over five random
initializations. The PINN is roughly five times
less accurate, but its errors remain at the $10^{-3}$
level,} which is more than adequate for most practical identification tasks.
The decisive difference, however, is in scalability. 
\revaa{This distinction becomes particularly important as the spatial dimension
increases.} \revzz{For a grid-based solver on a tensor-product mesh with $m$ points
per coordinate direction in $d$ space dimensions, the number of spatial degrees of
freedom is $N_x=m^{d}$, so the trajectory the adjoint stores and differentiates,
of size $N_x N_t$, grows exponentially in $d$.} \revaa{PINNs do not require a full
tensor-product discretization or storage of the discrete trajectory and therefore
avoid this particular source of exponential growth, although their cost may still
increase through the network capacity and collocation budget required in higher
dimensions. In the present 3D benchmark,} each adjoint iteration
sweeps forward and backward through the full 3D IMEX trajectory on a $129^3$
grid and caches the entire state history for the backward pass, so both its
runtime and its memory footprint grow with the spatial resolution; in four or
more dimensions, or on the fine grids needed to resolve thin interfaces, storing
and differentiating the full discrete trajectory becomes intractable. 

This is the central advantage of the PINN approach: on this 3D
benchmark it is already the cheaper option in wall-clock terms
(\revzz{$410\pm13$\,s versus $5825\pm350$\,s}) while remaining accurate to the
$10^{-3}$ level, and it extends to
regimes where assembling and storing a discretization-consistent adjoint is no
longer feasible.
Restarting the discrete adjoint from the PINN-recovered reaction network---the
``Adjoint restart'' history in Fig.~\ref{fig:ac3d_history_comparison}---begins
the backward sweeps already inside the basin of the true reaction, so the
\revzz{adjoint spends its budget improving rather than locating. It passes the
median accuracy the cold start finishes at after $2190\pm399$\,s, $38\%$ of the
cold cost, or $45\%$ counting the $410$\,s PINN that supplied the estimate.}

\subsection{Test 4: Scalar-viscosity identification for the 2D unsteady Navier--Stokes cylinder wake}
\label{subsec:cylinder_inverse}

The fourth benchmark treats the unsteady incompressible Navier--Stokes
equations past a circular cylinder at $\mathrm{Re}=100$, in the K\'arm\'an
vortex-shedding regime, with the unknown now the scalar kinematic viscosity
$\nu$ rather than a spatial field---the canonical problem of recovering a
single physical parameter from sparse wake observations, and the simplest
instance of the data-driven model calibration the paper ultimately targets.
Solving the cylinder wake by PINNs is itself a well-established inverse benchmark, treated
previously via dense wake velocity data~\cite{raissi2019physics},
vortex-induced vibration~\cite{raissi2019viv}, and sparse-measurement heat
transfer~\cite{cai2021heat}. Our setup is distinguished by combining an
explicitly resolved no-slip cylinder, sparse time-resolved probes, and a fully
unsteady regime in recovering the scalar viscosity that governs the shedding
dynamics.

The forward problem is set on the truncated fluid domain
$\Omega_f=([-3,10]\times[-3,3])\setminus\overline{B((0,0);0.5)}$ and governed by
\begin{equation}
\begin{aligned}
&\partial_t\mathbf{u}+(\mathbf{u}\!\cdot\!\nabla)\mathbf{u}+\nabla p-\nu\,\Delta\mathbf{u}=\mathbf{0},
& \mathbf{x}\in\Omega_f,\ t\in(0,T],\\
&\nabla\!\cdot\!\mathbf{u}=0, & \mathbf{x}\in\Omega_f,\ t\in(0,T],\\
&\mathbf{u}=(1,0),\ \mathbf{u}=\mathbf{0},\ v=0,
& \text{on inlet, cylinder, top/bottom walls}, \\
&\mathbf{u}(\mathbf{x},0)=\mathbf{u}_0(\mathbf{x}), & \mathbf{x}\in\Omega_f,
\end{aligned}
\label{eq:cylinder}
\end{equation}
with a natural outflow condition on $\{x=10\}$. In the abstract framework of
Section~\ref{sec:inverse} the unknown is the scalar $f\equiv\nu\in\mathbb{R}_{>0}$;
positivity is enforced by reparameterizing as $\theta_f=\log\nu$, so the
optimization variable is unconstrained and $\nu=\exp(\theta_f)$. The reference
value is $\nu^\star=10^{-2}$ (that corresponds to the Reynolds number $\mathrm{Re}=100$). The initial condition
$\mathbf{u}_0$ is not treated as unknown: a single warmup run generates the shedding state, which is then reused as the fixed initial
condition for every inverse run.

The forward problem is discretized by $\mathbb{P}_2$--$\mathbb{P}_1$ Taylor--Hood
finite elements on an unstructured triangular mesh with characteristic spacings
$h_{\mathrm{cyl}}=0.04$ near the cylinder, $h_{\mathrm{wake}}=0.08$ in the wake,
and $h_{\mathrm{far}}=0.5$ in the far field; time integration uses a second-order
IMEX SBDF2 scheme. Over the observation window $[0,T]$ with $T=5$, the solver is
marched with fixed step $\Delta t=5\times 10^{-3}$ ($N_t=1000$ steps); the discrete
adjoint differentiates this same trajectory
(~\ref{subsec:cylinder_implementation}).

The observation is sparse and time-resolved: both velocity components are sampled
at a $4\times 4$ Cartesian grid of $16$ probes covering
$[1,3]\times[-1,1]$ in the wake, recorded every $100$ time steps for
$N_{\mathrm{obs}}=10$ snapshots at $t\in\{0.5,1.0,\dots,5.0\}$, giving
$\mathbf{y}\in\mathbb{R}^{320}$ scalar observations (noise-free). The inverse problem~\eqref{eq:abstract_inverse} reduces to
\begin{equation}
\tilde{J}(\theta_f)
\;=\;\frac{1}{2}\sum_{k=1}^{N_{\mathrm{obs}}}\bigl\|P\,\mathbf{u}(\theta_f;t_k)-\mathbf{y}_k\bigr\|_2^{2},
\label{eq:cylinder_objective}
\end{equation}
with no explicit regularization ($\mathcal{R}\equiv 0$), where $P$ is the
$\mathbb{P}_2$ interpolation operator at the probe locations.

Two reconstructions are compared on this benchmark, both started from a
deliberately poor initial guess $\nu_0 = 10^{2}\nu^\star = 1.0$. The first is a
PINN utilizing a streamfunction--pressure network
$N_{\boldsymbol{\theta}_u}(\mathbf{x},t)$ outputting $(\psi,p)$ alongside a
trainable scalar $\theta_f = \log\nu$; the network is a multi-layer perceptron
(MLP) with three hidden layers of width $32$ and $\tanh$ activations. The
streamfunction--pressure parameterization enforces $\nabla\cdot\mathbf{u}=0$
exactly by construction ($\mathbf{u}=(\partial_y\psi,-\partial_x\psi)$). The
second is the discrete adjoint method, denoted ``Adjoint'', which optimizes
$\theta_f$ directly from the same scalar initialization with BFGS (SSBroyden's low-rank Hessian update is unnecessary for a single parameter).

\begin{figure}[!htbp]
    \centering
    \includegraphics[width=\textwidth]{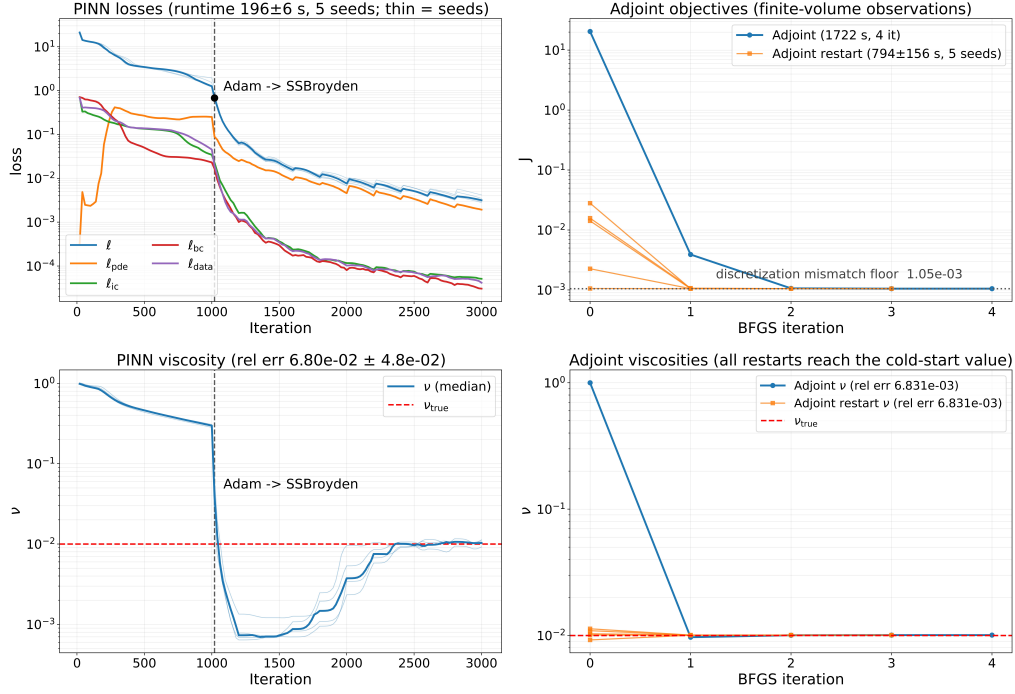}
    \caption{{\bf Test 4 - 2D unsteady Navier--Stokes cylinder wake.} Left: PINN loss
    and viscosity, Adam$\to$SSBroyden marked; right: adjoint objective and
    viscosity, cold-started and PINN-seeded. \revzz{Thin curves are five seeds, bold
    their median; every restart reaches the cold-started viscosity. The runtime in
    the legend is each restart's full wall-clock, of which the $567\pm133$\,s quoted
    in the text is the part needed to pass the accuracy the cold start finishes at.}}
    \label{fig:cylinder_history_comparison}
\end{figure}

\begin{figure}[!htbp]
    \centering
    \includegraphics[width=\textwidth]{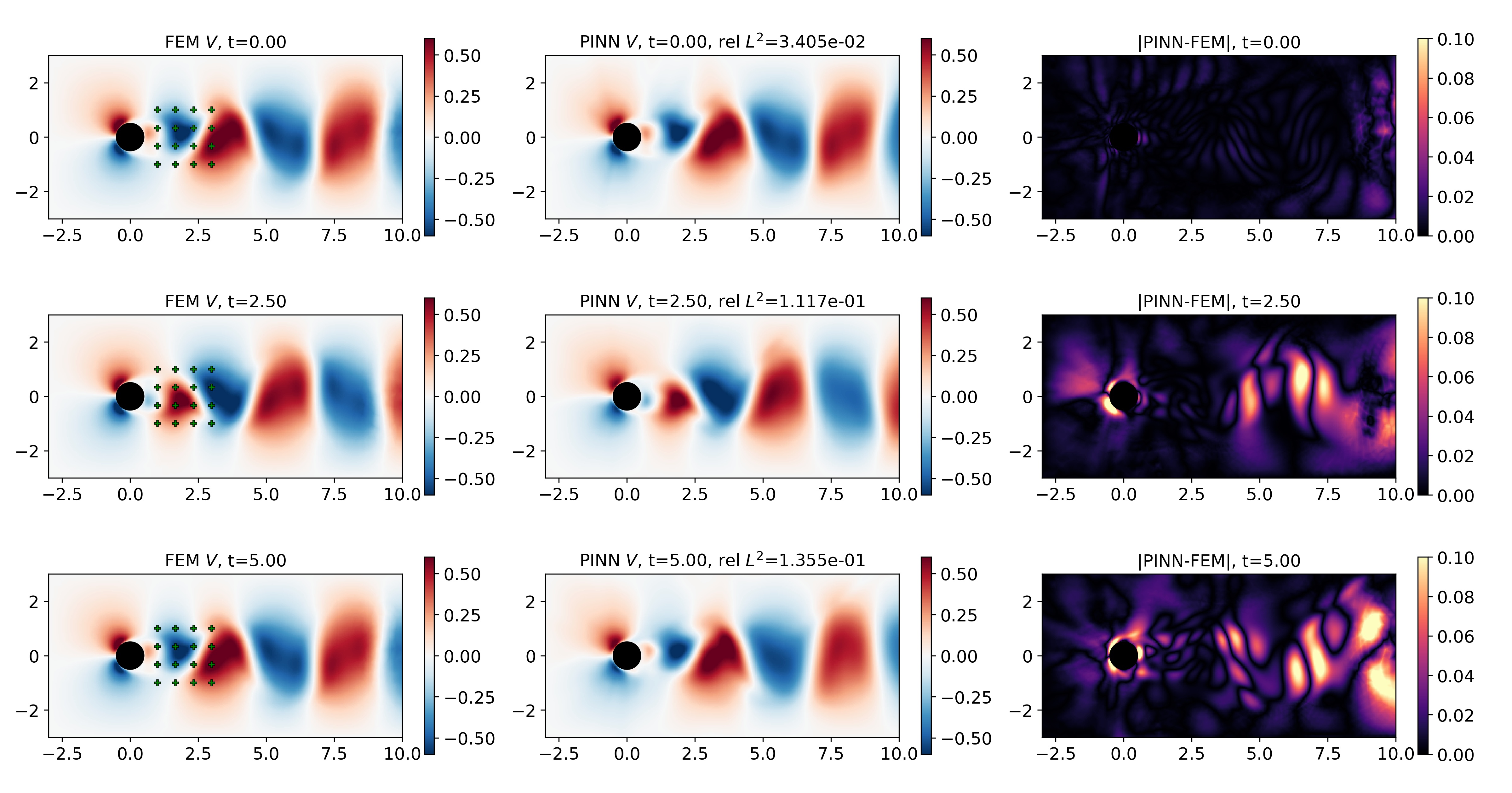}
    \caption{{\bf Test 4 - 2D unsteady Navier--Stokes cylinder
    wake.} Velocity-magnitude snapshots at $t\in\{0,\,2.5,\,5\}$. Left: FEM reference with the $16$
    wake probe locations (green crosses). Middle: inverse-PINN field at the
    converged viscosity. Right: pointwise absolute error. The probes constrain
    only a small wake window, yet drive the recovered viscosity to within
    \revzz{a few percent} relative error. \revzz{Single initialization.}}
    \label{fig:cylinder_velocity_snapshots}
\end{figure}

The convergence histories are reported in Fig.~\ref{fig:cylinder_history_comparison}.
\revzz{Both estimators are fitted to observations from an
independent, grid-converged solution. The baseline adjoint converges in $4$ BFGS
iterations, reaching $\varepsilon_\nu^{\mathrm{Adj}}=6.83\times 10^{-3}$ at a
wall-clock cost of $1722$\,s. Since the data come from a different discretization
than the one the inversion solves, the objective cannot fall below the mismatch
between the two: it flattens at $1.05\times10^{-3}$, the floor marked in
Fig.~\ref{fig:cylinder_history_comparison}, and the recovered viscosity absorbs
that mismatch.}
The cost is dominated by the $1000$-step IMEX-SBDF2 forward sweep that each BFGS evaluation
requires. \revzz{The PINN, by contrast, converges in $196\pm6$\,s---almost an
order of magnitude faster than the adjoint---to
$\varepsilon_\nu^{\mathrm{PINN}}=6.80\pm4.8\times10^{-2}$ over five random
initializations.} The viscosity trajectory in
the lower-left panel of Fig.~\ref{fig:cylinder_history_comparison} shows the
characteristic non-monotone behavior of a scalar PINN inversion: the SSBroyden
quasi-Newton phase first over-corrects to $\nu\approx 10^{-3}$ before recovering
to $\nu^\star$ once the auxiliary state network has been adequately fit by the
PDE-residual loss.
\revzz{A few percent error in a scalar recovered from $16$ sparse wake probes is
a satisfactory result for this benchmark, but the spread across initializations is
large --- the five runs place $\nu$ between $0.0092$ and $0.0113$ --- while the
state error of the trained network itself varies by only a few percent, so it is
the viscosity that is unstable, not the flow reconstruction.}

Restarting the discrete adjoint from a PINN-recovered viscosity collapses the
iteration count and the wall-clock cost without changing the answer.
\revzz{From each of the five PINN estimates, which span
$\varepsilon_\nu=0.9$--$13\%$, every restart returns the cold-started viscosity in $1$ to $3$ BFGS iterations,
passing the accuracy the cold start finishes at after $567\pm133$\,s; counting
the $196$\,s PINN that supplied the estimate, the two stages reach that accuracy
in $763$\,s against $1722$\,s, or $44\%$. The warm start therefore buys
cost rather than accuracy: this inverse problem has a single well-defined minimum,
and where the optimizer enters its basin changes only how quickly it arrives.}
For data-driven closure modeling of unsteady fluid models, where neither a one-shot
deterministic solve nor a soft-physics surrogate alone is fully satisfactory, this
hybrid is the practical recipe the paper recommends.

Figure~\ref{fig:cylinder_velocity_snapshots} shows the PINN's reconstructed
flow field over the full cylinder domain. The unsteady flow past a cylinder,
with the body and its no-slip surface resolved rather than reduced to a wake
patch, has remained a challenging problem~\cite{raissi2019physics,raissi2019viv,cai2021heat} for PINNs under sparse observations,
so we report the network's own reconstruction directly---the trained state
network $N_{\boldsymbol{\theta}_u}$ evaluated at the converged
viscosity---rather than only the re-simulated state error $\varepsilon_u$ of
Table~\ref{tab:inverse_summary}, which passes the recovered $\nu$ back through
the forward solver for a \revzz{controlled} cross-method comparison. The reconstruction is
demanding here because the only data available are the $16$ sparse wake
probes---unlike the other benchmarks, no full-field or terminal-state
observation constrains the network. The network nonetheless recovers the global
shedding structure from this limited information, with the pointwise error
concentrated in the near-wake, whereas in the other benchmarks, where richer
observations are available, the field is reconstructed near-perfectly.

\section{Conclusions}
\label{sec:conclusion}
We have compared adjoint-based optimization and physics-informed neural
networks for PDE-constrained inverse problems under a rigorous and well-controlled protocol. Both
methods were instantiated from a single abstract inverse-problem formulation,
with the domain, governing equation, and regularization held identical and the
optimizer, parameterization, and arithmetic precision matched wherever
applicable, \revzz{so that the measured differences are properties of the two
paradigms rather than of implementation choices}. This required care on both sides: the discretization
scheme was selected and linearized per benchmark for the adjoint, while the PINN
used a quasi-Newton optimizer, double precision, and residual resampling. The
comparison spans four benchmarks of increasing nonlinearity
and dimensionality---forcing identification for the unsteady 1D Burgers
equation, log-permeability identification for the steady 2D Darcy equation under
sparse noisy data, state-dependent reaction identification for the unsteady 3D
Allen--Cahn equation, and scalar-viscosity identification for the unsteady 2D
Navier--Stokes flow past a cylinder.

Across the benchmarks considered here, the representation of the unknown
emerges as a primary factor controlling the relative performance of the two
methods. When it is a state-independent forcing or a scalar parameter admitting a
low-dimensional, discretization-consistent representation, the adjoint method is
both the most accurate and the cheapest: it recovers the Burgers forcing near the
discretization floor, \revzz{identifies the cylinder viscosity to the accuracy its
observations support}, and solves the steady Darcy problem in seconds in a
prior-aligned space. Under sparse, noisy observations the three estimators become
indistinguishable, the residual gap set by the ill-posedness of the inverse problem
rather than the inversion algorithm, so the choice again reduces to efficiency and
favors the low-dimensional adjoint. The picture reverses once the unknown is
parameterized by a neural network, as required for state-dependent functionals
and data-driven closure: adjoint optimization converges poorly in the
high-dimensional weight space \revzz{and is prone to getting stuck in worse local minima}, whereas this representation is native to the PINN.

Time dependence and dimensionality add a second axis. \revaa{In particular, the discrete adjoint inherits the grid-induced curse of dimensionality of the underlying PDE discretization, whereas PINNs avoid the need to construct and store a full space-time grid.} Each adjoint iteration
stores and differentiates the entire discrete trajectory, so its runtime and
memory grow with the spatial resolution and the number of time steps and can become the dominant bottleneck in three-dimensional problems or over long time horizons. On the contrary, the PINN cost is
set by the network architecture and the collocation budget. On the 3D
Allen--Cahn problem the PINN is an order of magnitude cheaper and remains
accurate where a discretization-consistent adjoint is barely affordable. Where
\revzz{deterministic accuracy} is nonetheless required, the two methods combine: a
PINN warm start supplies a cheap, robust initial guess that the discrete adjoint
polishes \revzz{in a handful of iterations}. This
PINN-warm-started restart recovers the deterministic accuracy of the adjoint at
\revzz{$45\%$ and $44\%$ of the cold-started cost on the 3D Allen--Cahn reaction
and the unsteady cylinder viscosity respectively, counting the PINN that supplies
the estimate.}

These results give concrete guidance for method selection. Use the adjoint
method when a reliable discrete forward solver is available and the unknown is
low-dimensional---a scalar, a coefficient, or a coarse field. Use the PINN for
neural-field and state-dependent unknowns and for high-dimensional or
long-horizon time-dependent problems, where the adjoint's trajectory storage
becomes prohibitive. When both are feasible and high accuracy is required, use the
PINN-warm-started adjoint hybrid: the PINN supplies a low-cost approximation inside
the basin of attraction, and the discrete adjoint recovers high-fidelity accuracy
from it at a fraction of the cold-started cost.

A natural next step is to extend the
comparison to large-scale three-dimensional closure and calibration problems,
where the representation and cost trade-offs identified here are most
pronounced.

\section*{Declaration of competing interests}
The authors declare that they have no known competing financial interests or personal relationships that
could have appeared to influence the work reported in this paper.
\section*{Data Availability}
The source code, benchmark configurations, and data generated during this study are publicly available in the GitHub repository \url{https://github.com/zhangzhen117/adjoint_PINN_inverse_comparison}.

\section*{Acknowledgments}

This research was primarily supported by the Defense Advanced Research Projects Agency (DARPA) under the Automated Prediction Aided by Quantized Simulators (APAQuS) program, Grant No. HR00112490526. 
AA is a member of the GNCS group of INDAM. AA is partially supported by Progetto di Ateneo 2025: \emph{Bridging Analysis and Computation in Evolutionary PDEs.}
The authors would also like to thank Zhiwei Gao for implementing the EKI, and Shanqing Liu for valuable discussions.

\section*{Declaration of generative AI and AI-assisted technologies in the manuscript preparation process}
During the preparation of this work the authors used Claude in order to implement the computational code and polish the writing of the manuscript. After using this tool/service, the authors reviewed and edited the content as needed and take full responsibility for the content of the published article.
   
\bibliographystyle{elsarticle-num}
\bibliography{references}
\appendix
\section{Implementation Details}
\label{sec:appendix}
This appendix collects the discretization, parameterization, regularization, and optimization details deferred from Section~\ref{sec:benchmarks}. The common machinery shared by all four tests is described first in \ref{subsec:common_pinn}, followed by the test-specific configurations (\ref{subsec:burgers1d_implementation}--\ref{subsec:cylinder_implementation}). \revzz{The ablation studies are collected separately in~\ref{sec:ablation}.}
All computations use double-precision arithmetic throughout.

\subsection{Common PINN training pattern}
\label{subsec:common_pinn}

Unless stated otherwise, every PINN is trained in two stages. An Adam warmup with
a stepwise-decaying learning rate is run first for $\sim1000$ steps to bring the composite loss into the
basin of a good minimum, followed by a two-loop SSBroyden quasi-Newton phase. The outer
loop resamples the collocation and boundary points; the inner loop is on the resulting fixed batch, so the objective is deterministic within
each inner solve. Between outer iterations the inverse Hessian returned by SSBroyden is
symmetrized and accepted as the warm start for the next iteration only if it passes
a Cholesky positive-definiteness check, and is reset to the identity otherwise.
This restart-with-resampling pattern supplies curvature continuity across the
stochastic objective changes induced by resampling.

\begin{table}[ht]
\centering
\footnotesize
\setlength{\tabcolsep}{6pt}
\caption{PINN sampling budget per outer restart for the four benchmarks of
Section~\ref{sec:benchmarks}. The resampled column is randomly redrawn at every outer restart of the two-loop SSBroyden phase, while the fixed column remains constant across training.}
\label{tab:pinn_sampling}
\begin{tabular}{clp{4.0cm}l}
\toprule
Test & Benchmark & Resampled & Fixed \\
\midrule
1 & 1D Burgers       & 20{,}000 PDE + 512 BC                   & 512 IC + 512 data            \\
2 & 2D Darcy         & 2000 PDE + 200 BC                       & 4225 data                    \\
3 & 3D Allen--Cahn   & 20{,}000 PDE + 5000 BC \newline + 80{,}000 data + 80{,}000 IC & --- \\
4 & 2D NS cylinder   & 10{,}000 PDE + 6000 BC                  & 24{,}745 IC + 320 data \\
\bottomrule
\end{tabular}
\end{table}

All PINN loss terms carry unit weight unless noted, except the regularization
term, whose weight is the same $\gamma$ used by the adjoint method; in every case the state field and the unknown
parameter are represented by separate networks. The per-benchmark sampling
budget is summarized in Table~\ref{tab:pinn_sampling}, where the
\emph{resampled} columns list points redrawn uniformly at random at every outer
quasi-Newton restart and the \emph{fixed} columns list points held constant
throughout training. For the 3D Allen--Cahn benchmark the terminal-state
observation is the full $129^3$ field; rather than enforcing the data
misfit on all $\sim2.1\times10^6$ nodes at once, a random subset is resampled
each outer loop, which both bounds the per-iteration cost and acts as a
stochastic regularizer against overfitting to any fixed point cloud.

\subsection{Test 1: 1D viscous Burgers}
\label{subsec:burgers1d_implementation}
\paragraph{Forward solver}
The forward problem~\eqref{eq:burgers1d} is discretized on a uniform periodic grid
of $N_x=512$ nodes with spacing $\Delta x=L/N_x$. The first derivative and the
Laplacian use centered second-order stencils assembled as periodic circulant
matrices $D,L_x\in\mathbb{R}^{N_x\times N_x}$. With
$F(u;f)=-u\odot(Du)+\nu L_x u+f$, the Crank--Nicolson update at step $n$ is enforced
as the residual
\begin{equation}
R^n(\mathbf{u}^{n-1},\mathbf{u}^{n};f)
=\mathbf{u}^{n}-\mathbf{u}^{n-1}
-\tfrac{\Delta t}{2}\bigl[F(\mathbf{u}^{n-1};f)+F(\mathbf{u}^{n};f)\bigr]=0,
\qquad n=1,\dots,N_t,
\label{eq:burgers_cn_residual}
\end{equation}
with $\Delta t=10^{-2}$ ($N_t=100$ steps to $T=1$). Each implicit step is solved for
$\mathbf{u}^{n}$ by Newton's method with the analytic Jacobian
$J_F(u)=-\mathrm{diag}(Du)-\mathrm{diag}(u)\,D+\nu L_x$ to an $\ell_\infty$ residual
tolerance of $10^{-10}$ (at most $20$ iterations, explicit-Euler predictor). A
mesh-refinement study confirms second-order accuracy in both $\Delta x$ and
$\Delta t$, giving the discretization floor
$\varepsilon_u^{\mathrm{solver}}=4.19\times 10^{-5}$ at the production resolution.
\paragraph{Discrete adjoint}
Differentiating the residual~\eqref{eq:burgers_cn_residual} gives the step Jacobians
$A^n=\partial R^n/\partial\mathbf{u}^{n-1}=-I-\tfrac{\Delta t}{2}J_F(\mathbf{u}^{n-1})$ and
$B^n=\partial R^n/\partial\mathbf{u}^{n}=I-\tfrac{\Delta t}{2}J_F(\mathbf{u}^{n})$.
The adjoint variables solve the backward sweep
$(B^{N_t})^\top\boldsymbol{\lambda}^{N_t}=-\Delta x\,(\mathbf{u}^{N_t}-\mathbf{y})$
and $(B^{n})^\top\boldsymbol{\lambda}^{n}=-(A^{n+1})^\top\boldsymbol{\lambda}^{n+1}$
for $n=N_t-1,\dots,1$, and since $\partial R^n/\partial f=-\Delta t\,I$, the gradient
with respect to the grid forcing is
$\partial J/\partial f=-\Delta t\sum_{n=1}^{N_t}\boldsymbol{\lambda}^{n}$. The
reduced gradient is verified against centered finite differences of the discrete
objective along randomly sampled perturbation directions, agreeing to a maximum
relative error of $\approx 10^{-8}$.

\subsection{Test 2: 2D Darcy flow}
\label{subsec:darcy_implementation}

\paragraph{Forward solver}
The Darcy problem~\eqref{eq:darcy} is discretized by continuous
$\mathbb{P}_2$ Lagrange finite elements on a structured triangular mesh obtained by
splitting each cell of a $32\times 32$ grid into two triangles, giving
$N_x=4225$ degrees of freedom. Element integrals use a $7$-point Gauss rule of order
$5$. The permeability enters the stiffness matrix elementwise as
$k_e=\exp(f_e)$ via a precomputed unit-permeability assembly that is reweighted at
each solve; homogeneous Dirichlet conditions are imposed by interior-degree-of-freedom
restriction and the linear system is solved by a sparse direct factorization. A
manufactured-solution study confirms the expected $\mathbb{P}_2$ rates ($O(h^3)$ in
$L^2$), giving $\varepsilon_u^{\mathrm{solver}}=7.45\times 10^{-6}$ at the production
mesh.

\paragraph{Prior and parameterization}
The reference log-permeability and the priors of the adjoint and
EKI~\cite{iglesias2013ensemble} methods are samples from a Gaussian random field
expressed through a truncated Karhunen--Lo\`eve expansion,
\begin{equation}
f(\mathbf{x})=\sum_{i=1}^{r}\xi_i\,\sqrt{\lambda_i}\,\varphi_i(\mathbf{x}),
\qquad \xi_i\sim\mathcal{N}(0,1),
\end{equation}
with tensor-product cosine eigenfunctions
$\varphi_i(\mathbf{x})=\cos(\pi k_{i,1}x_1)\cos(\pi k_{i,2}x_2)$ indexed by integer
wavenumbers $k_i=(k_{i,1},k_{i,2})\in\mathbb{Z}_{\ge 0}^{\,2}$, and corresponding
eigenvalues
$\sqrt{\lambda_i}=\sigma\bigl(\pi^{2}\|k_i\|^{2}+\tau^{2}\bigr)^{-\alpha/2}$ of a
Mat\'ern-type covariance with production parameters $\tau=3$, $\alpha=2$, and
amplitude $\sigma=1$. We truncate the expansion at $r=128$ terms (ordered by
descending $\lambda_i$) to obtain a finite-dimensional approximation of the
infinite-dimensional field. Let
$\Phi_{\mathrm{sc}}\in\mathbb{R}^{r\times N_x}$ be the scaled eigenfunction matrix
with rows $\sqrt{\lambda_i}\,\varphi_i$ evaluated at the $N_x$ FEM degrees of
freedom, so the discrete log-permeability is
$\mathbf{f}=\Phi_{\mathrm{sc}}^{\top}\boldsymbol{\xi}\in\mathbb{R}^{N_x}$ and the
chain rule gives $\nabla_{\boldsymbol{\xi}}J=\Phi_{\mathrm{sc}}\nabla_{\mathbf{f}}J$.
The adjoint and EKI methods both optimize the $r=128$ KL coefficients
$\boldsymbol{\xi}$ directly. The PINN instead represents the log-permeability by a
neural field $N_{\boldsymbol{\theta}_f}(\mathbf{x})$ alongside the state network
$N_{\boldsymbol{\theta}_u}(\mathbf{x})$.

\paragraph{Discrete adjoint}
For this steady problem the adjoint reduces to a single linear solve with the
transpose of the converged-state Jacobian, and the reduced gradient follows by
the chain rule through the KL map, $\nabla_{\boldsymbol{\xi}}J
=\Phi_{\mathrm{sc}}\nabla_{\mathbf{f}}J$. The gradient is verified against
centered finite differences of the discrete objective along randomly sampled
perturbation directions, agreeing to a maximum relative error of
$\approx 10^{-7}$.

\paragraph{Ensemble Kalman inversion}
To implement EKI, we write the observation model in the equivalent form
\begin{equation}
    y = \mathcal G(\boldsymbol{\xi}^{\star})+\eta,
    \qquad
    \eta\sim\mathcal N(0,\Gamma),
    \qquad
    \Gamma=\sigma^2 I,
\end{equation}
where the unknown field is parametrized by the random coefficients
\(\boldsymbol{\xi}\) through the expansion \(f=f(\boldsymbol{\xi})\), and
\begin{equation}
    \mathcal G(\boldsymbol{\xi})
    :=
    \mathcal O\bigl(u(f(\boldsymbol{\xi}))\bigr).
\end{equation}
Thus, \(\mathcal G\) includes the random-field expansion, the PDE/FEM solve, and
the sparse observation operator, and \(\sigma\) is the observation noise level
used in the experiment.

EKI approximates the inverse problem by evolving an ensemble
\(\{\boldsymbol{\xi}_n^{(j)}\}_{j=1}^{J}\), where \(J\) is the ensemble size and
\(n\) is the iteration index. The initial ensemble is sampled from the prior,
\begin{equation}
    \boldsymbol{\xi}_0^{(j)}\sim \mathcal{N}(0, I),
    \qquad j=1,\ldots,J.
\end{equation}
At iteration \(n\), the forward map is evaluated for each ensemble member,
\(g_n^{(j)}=\mathcal G(\boldsymbol{\xi}_n^{(j)})\), and the ensemble means are
\begin{equation}
    \bar{\boldsymbol{\xi}}_n
    =\frac{1}{J}\sum_{j=1}^{J}\boldsymbol{\xi}_n^{(j)},
    \qquad
    \bar g_n
    =\frac{1}{J}\sum_{j=1}^{J}g_n^{(j)}.
\end{equation}
The empirical cross- and observation covariances are
\begin{equation}
    C_n^{\xi g}
    =\frac{1}{J-1}\sum_{j=1}^{J}
    \bigl(\boldsymbol{\xi}_n^{(j)}-\bar{\boldsymbol{\xi}}_n\bigr)
    \bigl(g_n^{(j)}-\bar g_n\bigr)^\top,
    \qquad
    C_n^{gg}
    =\frac{1}{J-1}\sum_{j=1}^{J}
    \bigl(g_n^{(j)}-\bar g_n\bigr)
    \bigl(g_n^{(j)}-\bar g_n\bigr)^\top,
\end{equation}
and the ensemble is updated by
\begin{equation}
    \boldsymbol{\xi}_{n+1}^{(j)}
    =\boldsymbol{\xi}_{n}^{(j)}
    +C_n^{\xi g}\bigl(C_n^{gg}+\Gamma\bigr)^{-1}
    \bigl(y_n^{(j)}-g_n^{(j)}\bigr),
    \qquad j=1,\ldots,J,
\end{equation}
with perturbed observations \(y_n^{(j)}=y+\eta_n^{(j)}\),
\(\eta_n^{(j)}\sim\mathcal N(0,\Gamma)\); for the deterministic variant we set
\(y_n^{(j)}=y\) for all members. After \(N\) iterations the reconstructed
coefficients are the ensemble mean
\(\widehat{\boldsymbol{\xi}}=\tfrac{1}{J}\sum_{j=1}^{J}\boldsymbol{\xi}_N^{(j)}\),
and the reconstructed field is \(\widehat f = f(\widehat{\boldsymbol{\xi}})\).

\paragraph{Regularization}
The adjoint regularizer is a discrete $H^1$ smoothness penalty applied on the
piecewise-constant log-permeability field over the FEM triangulation,
\begin{equation}
\mathcal{R}(\mathbf{f})=\frac{\gamma}{2}\!\!\sum_{(e_i,e_j)\in\mathcal{E}_{\mathrm{int}}}\!\!
\frac{(f_{e_i}-f_{e_j})^2}{\ell_{ij}},\qquad \gamma=\gamma^{\ast}(\sigma),
\end{equation}
\revzz{where the weight is selected once per noise level from the grid of
Table~\ref{tab:darcy_noise} and shared by both gradient-based methods;
$\gamma^{\ast}=10^{-2}$ at the $\sigma=1\%$ of the production runs.}
Here $\mathcal{E}_{\mathrm{int}}$ is the set of pairs of adjacent triangles sharing
an interior edge, $\ell_{ij}$ is the length of the shared edge, and $f_{e_i}$ is the
log-permeability value on triangle $e_i$. Equivalently, this is the quadratic form
$\tfrac{\gamma}{2}\mathbf{f}^{\top}L_{\mathrm{graph}}\mathbf{f}$ associated with the
weighted graph Laplacian of the dual mesh; on piecewise-constant fields it is a
finite-volume discretization of the continuous $H^1$ seminorm
$\int_{\Omega}|\nabla f|^{2}\,d\mathbf{x}$.

\subsection{Test 3: 3D Allen--Cahn}
\label{subsec:ac3d_implementation}

\paragraph{Forward solver}
The 3D Allen--Cahn problem~\eqref{eq:ac3d} is discretized on a uniform
tensor-product grid of $N_x=N_y=N_z=129$ nodes by the seven-point finite-difference
Laplacian (a Kronecker sum of one-dimensional stencils) with ghost-point Neumann
closure on all faces. Time integration is a second-order IMEX scheme: a first-order
IMEX Euler step initializes the multistep method, after which the uniform steps use
IMEX BDF2/AB2,
\begin{equation}
\bigl(3I-2\Delta t\,\varepsilon^2 L_h\bigr)\mathbf{u}^{n+1}
=4\mathbf{u}^{n}-\mathbf{u}^{n-1}+2\Delta t\bigl(2f(\mathbf{u}^{n})-f(\mathbf{u}^{n-1})\bigr),
\end{equation}
with $\Delta t=5\times 10^{-2}$. Because the implicit operator is the Neumann
Laplacian, each implicit solve is diagonalized by discrete cosine transforms (DCT) with
eigenvalues $\lambda_k=-\tfrac{4}{h^2}\sin^2\!\bigl(\tfrac{\pi k}{2(N-1)}\bigr)$
summed across directions, so each step is a transform--divide--inverse-transform
operation. The discretization floor at this resolution is
$\varepsilon_u^{\mathrm{solver}}=4.40\times 10^{-4}$.

\paragraph{Discrete adjoint}
The adjoint uses the same IMEX structure with $f=N_{\boldsymbol{\theta}_f}$, caching
the current and previous reaction evaluations, the derivative $\partial f/\partial u$
at each visited state, and the step type. The backward sweep starts from the terminal
mismatch $\boldsymbol{\lambda}^{N_t}=\Delta x\Delta y\Delta z\,(\mathbf{u}^{N_t}-\mathbf{y})$ and
applies the same DCT-diagonalized solver to adjoint loads; parameter gradients are
accumulated by backpropagation through the reaction network at the stored states. The
adjoint gradient is validated by random-direction finite differences (maximum relative
error $\approx 1.2\times 10^{-4}$).

\subsection{Test 4: 2D unsteady Navier--Stokes cylinder wake}
\label{subsec:cylinder_implementation}

\paragraph{Forward solver}
The cylinder problem~\eqref{eq:cylinder} is discretized by
$\mathbb{P}_2$--$\mathbb{P}_1$ Taylor--Hood finite elements on an unstructured
triangular mesh of the truncated fluid domain
$\Omega_f=([-3,10]\times[-3,3])\setminus\overline{B(0,0.5)}$, with characteristic
spacings $h_{\mathrm{cyl}}=0.04$ near the cylinder, $h_{\mathrm{wake}}=0.08$ in
the wake, and $h_{\mathrm{far}}=0.5$ in the far field. Element integrals use a
$7$-point Gauss rule of order $5$. The inlet and cylinder Dirichlet conditions
are imposed by interior-degree-of-freedom restriction; the top and bottom walls
carry a slip-style $v=0$ constraint; the outflow at $\{x=10\}$ is the natural
condition of the weak form (no constraint). The pressure null space is removed
by a small mass-matrix regularization with parameter $\epsilon=10^{-10}$. Time
integration is a second-order IMEX SBDF2 scheme: convection is extrapolated by
SBDF2 from the two previous steps while diffusion is treated implicitly, so each
time step requires one sparse direct factorization of the velocity stiffness
block.

The reference data are generated in two stages. A warmup integration to
$T_{\mathrm{warm}}=80$ uses adaptive time-stepping with CFL target $0.9$ and a symmetry-breaking blob perturbation of amplitude
$0.1$ injected at $t=8$, advancing the flow into the K\'arm\'an
shedding cycle. The final state is then frozen as $\mathbf{u}_0$ and the
observation window $[0,T]$ with $T=5$ is marched with fixed step
$\Delta t=5\times 10^{-3}$ ($N_t=1000$ steps). The same saved $\mathbf{u}_0$ is
reused as the initial condition for every inverse run, so $\mathbf{u}_0$ is not
treated as an unknown.

\paragraph{Discrete adjoint}
The discrete adjoint linearizes the converged forward trajectory and is
implemented as a backward sweep through the cached IMEX-SBDF2 states. The
optimization variable is the scalar $\theta_f=\log\nu$, so $\nu=\exp(\theta_f)$
and viscosity positivity is enforced by construction; the chain rule gives
$dJ/d\theta_f=\nu\,dJ/d\nu$. The gradient $dJ/d\nu$ is the sum, across the $1000$
time steps and $10$ observation snapshots, of the inner products of the cached
adjoint velocity multipliers with the diffusion operator applied to the
state---a single scalar derivative per step. The adjoint is verified against
centered finite differences with $\delta=10^{-6}$ at viscosity offsets
$\nu/\nu^\star\in\{0.7,0.85,1.2,1.5\}$; the worst-case
relative error is below $10^{-7}$.

\revzz{\section{Ablation Studies}
\label{sec:ablation}}

\revzz{This appendix reports two sensitivity studies: the optimizer, on the simplest benchmark, and PINN's architecture, on the unsteady three-dimensional benchmark where the PINN's advantage is largest.}

\revzz{\subsection{PINN optimizer comparison}
\label{subsec:burgers_optimizer}}

\revzz{The Burgers PINN was retrained with SOAP~\cite{vyas2024soap} and with Adam
in place of SSBroyden, holding the network, collocation budget, data and precision fixed.
Each first-order run is driven to saturation rather than to a fixed budget: the
learning rate quoted below is an initial value, tuned per optimizer over
$\{10^{-3},3{\times}10^{-3},10^{-2},3{\times}10^{-2}\}$, decayed on plateau by a
factor $0.3$ after $5000$ steps without improvement and terminated once it falls to
$5\times10^{-6}$ of that value. The tuned rate for Adam, $10^{-2}$, is therefore
larger than the value usually quoted for PINN training because it starts a decaying
schedule rather than fixing a step size. The
SSBroyden row is the production run of Section~\ref{sec:benchmarks}.}

\revzz{
\begin{table}[!htbp]
\centering
\color{zzrev}
\footnotesize
\caption{\revzz{PINN optimizer comparison on the Burgers benchmark.
Mean\,$\pm$\,standard deviation over five random initializations, each run to
saturation at its own tuned learning rate.}}
\label{tab:burgers_optimizer}
\begin{tabular}{lcccc}
\toprule
Optimizer & lr & $\varepsilon_f$ & $t$ [s] & steps \\
\midrule
SSBroyden & ---                & $\mathbf{1.04{\pm}0.06{\times}10^{-4}}$ & $\mathbf{195{\pm}4}$ & $\sim\!4{\times}10^{3}$ \\
SOAP      & $3{\times}10^{-3}$ & $7.56{\pm}2.2{\times}10^{-4}$           & $1424{\pm}92$        & $1.8{\times}10^{5}$ \\
Adam      & $10^{-2}$          & $5.33{\pm}1.4{\times}10^{-3}$           & $1803{\pm}358$       & $3.3{\times}10^{5}$ \\
\bottomrule
\end{tabular}
\end{table}
}

\revzz{SSBroyden is $7$ times more accurate than SOAP and $51$ times more accurate
than Adam, at one seventh and one ninth of the cost respectively; the first-order
methods need on the order of $10^{5}$ updates to saturate against some
$4\times10^{3}$ quasi-Newton iterations.}

\revzz{SSBroyden maintains a dense $N_f\times N_f$ inverse-Hessian approximation,
so its storage is $8N_f^{2}$ bytes and each update is $\mathcal{O}(N_f^{2})$:
$10.6$\,MB at the $N_f=1153$ used here, $800$\,MB at $10^{4}$ parameters and
$80$\,GB at $10^{5}$. It is negligible at the network sizes of this study but
would be prohibitive for substantially larger networks, where a limited-memory
method would be required.}

\revzz{\subsection{Allen--Cahn: sensitivity to the PINN architecture}
\label{subsec:ac_arch}}

\revzz{Three axes are varied around the production setting, one at a time: the
size of both networks together, the activation --- including the layer-wise
locally adaptive form of~\cite{jagtap2020locally} --- and a random Fourier
embedding~\cite{tancik2020fourier,wang2021eigenvector} of the state network,
applied there only because the unknown $f(u)$ is a scalar function of one
variable. ReLU-type activations are excluded throughout, since the
residual requires $\partial_{xx}u$ and their second derivative vanishes.}

\revzz{All configurations are trained under one stopping rule --- outer epochs
continue until the loss improves by less than $10^{-4}$ over five consecutive
epochs, to a ceiling of $80$ --- rather than to the $15$ epochs of the production
runs, at which the production configuration is close to flat while the slower
variants still improve by $10$--$20\%$ per epoch. Three random initializations per
configuration.}

\revzz{
\begin{table}[!htbp]
\centering
\color{zzrev}
\footnotesize
\setlength{\tabcolsep}{4pt}
\caption{\revzz{Allen--Cahn PINN architecture study. Relative error of the recovered
reaction, mean\,$\pm$\,standard deviation over three random initializations, all
trained under the same stopping rule. ``Epochs'' is the number actually used, of a
ceiling of $80$; configurations that reach the ceiling were still improving when
stopped and their errors are upper bounds.}}
\label{tab:ac_arch}
\begin{tabular}{llccrr}
\toprule
Axis & Configuration & $N_{\mathrm{params}}$ & $\varepsilon_f$ & epochs & $t$ [s] \\
\midrule
          & production, $32{\times}3$ tanh & $3458$  & $\mathbf{1.27{\pm}0.29{\times}10^{-3}}$ & $80$ & $2204$ \\
\midrule
size      & small, $16{\times}2$           & $690$   & $1.05{\pm}0.35{\times}10^{-1}$ & $24$ & $\mathbf{171}$ \\
          & large, $64{\times}4$           & $17218$ & $1.34{\pm}0.10{\times}10^{-3}$ & $70$ & $26832$ \\
\midrule
activation& silu                           & $3458$  & $7.18{\pm}7.4{\times}10^{-3}$  & $80$ & $2272$ \\
          & adaptive (L-LAAF)              & $3463$  & $1.54{\pm}2.0{\times}10^{-1}$  & $52$ & $1389$ \\
\midrule
Fourier   & $m=8$, scale $1$               & $3842$  & $5.94{\pm}0.95{\times}10^{-3}$ & $78$ & $2367$ \\
          & $m=8$, scale $5$               & $3842$  & $5.02{\pm}5.0{\times}10^{1}$   & $35$ & $1126$ \\
\bottomrule
\end{tabular}
\end{table}
}

\revzz{The production configuration is the best of the seven. Increasing capacity
five-fold changes the error by $5\%$, within the seed spread, while costing twelve
times the wall-clock; reducing it to a fifth costs a factor of $83$. The chosen
size therefore sits at the point where accuracy saturates.}

\revzz{Replacing tanh by silu costs a factor of $5.7$ and the adaptive activation two
orders of magnitude, while the Fourier embedding costs $4.7$ at unit scale and
fails outright at scale $5$: differentiating $\sin(Bx)$ twice brings down a factor
$B^2$, amplifying precisely the curvature terms that the second-order residual
already makes stiff. On a second-order problem with a smooth solution, spectral
enrichment of the input is a liability rather than a benefit.
The wider reading of this table is that the PINN's accuracy does depend on
its architecture, and that arriving at a good one requires a search of the kind
carried out here. This is not an asymmetry between the two paradigms but the same
requirement in a different guise: the adjoint results rest equally on choices that
had to be designed rather than assumed --- the discretization scheme, the
resolution, and problem-specific structure such as the discrete cosine
diagonalization that makes the Allen--Cahn implicit solve matrix-free --- so what
differs is where the design effort is spent, not whether it is spent.}

\end{document}